\documentclass[twocolumn,showkeys,preprintnumbers,amsmath,amssymb]{revtex4}
\usepackage{bm}
\usepackage{dcolumn}
\usepackage{epsfig}
\usepackage{amsfonts}
\usepackage{color}
\usepackage{natbib}
\bibpunct{}{}{,}{s}{}{}

\begin{document}

\title{\raggedright \Large {\bf Equation-Free Particle-Based Computations: 
Coarse Projective Integration and Coarse Dynamic Renormalization in 2D}}

%
%

\author{\bf \hspace{-3.55cm} Yu Zou}
\thanks{Department of Chemical Engineering and PACM, Princeton University.}
\author{\bf Ioannis G. Kevrekidis}%
 \thanks{To whom correspondence should be addressed. Tel.: 1-609-258-2818. E-mail: yannis@princeton.edu. Department of Chemical Engineering and PACM, Princeton University.}
\author{\bf Roger G. Ghanem}
\thanks{Department of Civil Engineering, The University of Southern 
California.}
\affiliation{\hspace{-0.65cm} Department of Chemical Engineering and PACM, Princeton University, Princeton, NJ 08544
\\
Department of Civil Engineering, The University of Southern 
California, Los Angeles, CA 90089}%

\begin{abstract}

Equation-free approaches have been proposed in recent years 
for the computational study of multiscale phenomena  
in engineering problems where evolution equations for the coarse-grained,
system-level behavior are not explicitly available. 
In this paper we study the dynamics of a
diffusive particle system in a laminar shear flow,
described by a two-dimensional Brownian motion; in particular,
we perform coarse projective integration and demonstrate 
the particle-based computation of coarse self-similar 
and asymptotically self-similar solutions for this problem.
We use marginal and conditional Inverse Cumulative Distribution Functions
(ICDFs) as the macroscopic observables of the evolving particle distribution.
\end{abstract}


\maketitle
\noindent{\bf 1. Introduction}
\vspace{0.2cm}

Multiscale phenomena arise naturally in science and engineering. 
The ability to properly resolve such phenomena and propagate their influence 
across scales underpins the predictive value of mathematical and physics-based models. 
In the case of multiscale systems for whose macroscopic behavior no explicit 
coarse-grained, macroscopic equations are available, 
a computer-assisted approach, referred to as the {\it Equation-Free Framework}
\cite{Kevrekidis:03,Kevrekidis:04} has been recently proposed.
Equation-free methods numerically evolve the coarse-scale 
behavior through appropriately designed short computational
experiments performed by the fine-scale (microscopic, stochastic,
agent-based) models.
In this paper we will demonstrate the use of two such methods: Coarse
Projective Integration (CPI) and Coarse Dynamic Renormalization (CDR).

Particle-based simulators are the fine-scale description of choice
for a variety of problems exhibiting multiscale behavior; such problems
range from Stokesian and Brownian dynamics to the Monte Carlo modeling of
microorganism locomotion, the mixing of passive scalars 
by turbulent velocity fields and even particle filtering applications.
The purpose of Coarse Projective Integration is to accelerate
the computational evolution of coarse-grained observables of microscopic
simulators; it has been successfully used in the past to accelerate 
computations of the collective evolution of {\it spatially one-dimensional}
random particle distributions \cite{Gear:01,Setayeshgar:03}.
In those examples, the first motivated by hydrodynamics and the second by 
bacterial chemotaxis, the coarse-scale observable was the cumulative distribution 
function (CDF) of the particle
positions, synthesized from snapshots of the fine-scale simulation.
For one-dimensional problems in space,
the functional inverse of the CDF (ICDF) is projected onto a suitable basis set 
consisting of orthonormal polynomials \cite{Gear:01} or POD (Proper Orthogonal Decomposition) modes (obtained through Singular Value Decomposition (SVD) 
of the sorted particle positions \cite{Setayeshgar:03}).  
The particle-level model is used to evolve the ICDF, and short time series of the
coefficients of its projection on the appropriate basis are thus collected;
these short time series are used to estimate the time derivatives of the
coefficient evolution.
These estimates are then used in the context of traditional continuum integration 
algorithms (such as Euler \cite{Gear:02}, or Adams-Bashforth \cite{Faires:93}) to ``project" the
coefficients forward in time, into the future (thus the term ``projective integration").
To repeat the procedure, a ``lifting" step is required: fine-scale states, i.e., particle 
positions whose ICDF is consistent with the projected coefficient values, are generated.
Since detailed microscopic evolution over the duration of the projective step
has been avoided, the procedure has the potential to alleviate the burden of
full fine-scale simulation.
Preliminary discussions on the stability and accuracy of these schemes can be found in 
\cite{GearA:01,Gear:02}.

Beyond CPI, equation-free computational protocols can also be used 
in coarse-grained fixed point and bifurcation analyses, to compute stationary
states of the coarse-grained system dynamics and their parametric dependence.
In these analyses, the action of operators on the coarse-grained observables
is deduced from appropriately initialized computational experiments
with the fine scale models. 
Matrix-free implementations of contraction mappings, like Newton's method, have been used to compute fixed points of unavailable coarse-scale models
for kinetic Monte Carlo (e.g., \cite{Makeev:02}, epidemiology),
Brownian dynamics \cite{SiettosA:03} or Lattice Boltzmann \cite{Gear:02,Theodoropoulos:04}
fine-scale simulators.
The same approach can be used to evolve 
effective medium (homogenized) descriptions of reaction-transport problems
\cite{Runborg:02,Xiu:04} based on {\it short bursts} of finely resolved simulation.

For multiscale systems exhibiting {\it scale invariance} at the macroscopic level,
renormalization techniques can be used to solve for self-similar solutions
and their scalings \cite{Barenblatt:96}, see also \cite{Brandt:01,Chorin:03}.
Recently, {\it dynamic} renormalization
(e.g., \cite{Mclaughlin:86,LemesurierA:88,LemesurierB:88})
was used in conjuction with equation-free computation to obtain 
coarse-grained self-similar solutions
using short bursts of fine-level, direct simulation \cite{Chen:03}. 
This {\it coarse dynamic renormalization} (CDR) method
finds macroscopically self-similar solutions with the help of
{\it template functions}  
\cite{Rowley:00,Aronson:01,Siettos:03,Rowley:03}.
Observing the macroscopic solutions in a co-expanding 
(or co-collapsing) frame of reference, we seek 
{\it steady states} in the new frame; fixed point equation-free
algorithms can be used for this task.

Studying the evolution of particle distributions using their ICDF as
an observable is convenient for one-dimensional problems in space,
where suitable bases for representing one-dimensional monotonic
curves over finite one-dimensional domains are 
readily available (see the examples in \cite{Gear:01,Setayeshgar:03,Chen:03}). 
The extension to corresponding observables in more than
one dimension, however, is nontrivial: 
the CDF itself, not being a bijective mapping, does not have an inverse. 
In this case, operations on the CDF may be implemented by identifying 
a suitable set of basis functions for two or three-dimensional CDFs.
This requires finding multidimensional orthonormal
polynomial approximations for monotonic bounded 
functions with infinite support. 
%
%
In this paper we use an alternative approach, representing a multidimensional CDF 
in terms of its marginal and one-dimensional conditional distributions.
In this manner, multidimensional problems are converted into a collection of
one-dimensional problems the solution of which can be obtained using standard approaches.
Equation-free algorithms such as CPI and CDR can thus readily be 
extended to problems involving multidimensional coarse-grained observables.  
Preliminary results for the self-similar case have been reported elsewhere \cite{Zou:05}.

The paper is organized as follows:
In Section 2 a coarse time-stepper is constructed in terms of the marginal and 
conditional ICDFs in multidimensional particle systems.  
Use of this time-stepper in CPI and CDR is formulated 
in Sections 3 and 4, respectively. 
Our illustrative Brownian particle system in a Couette flow is described
in Section 5, and its analytical
self-similar and asymptotically self-similar solutions are presented.
These two cases are  
then used in Section 6 to illustrate equation-free computational
procedures, and the direct, particle-based computational results are compared with 
the analytical solutions.  
We conclude with a brief summary and discussion in Section 7.

\vspace{0.5cm}
\noindent{\bf 2. A Coarse Time-Stepper for Multidimensional Random Particle Systems} 
\vspace{0.2cm}

\begin{figure}
\centerline{\epsfig{figure=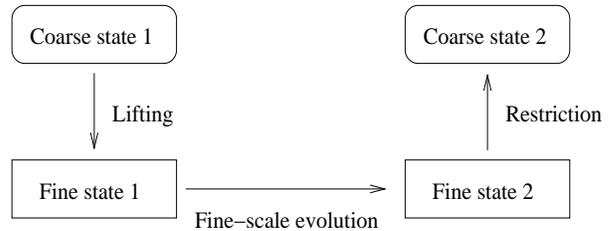,height=30mm,width=80mm,angle=0}}
\caption{A schematic of the coarse time-stepper.}
\label{cts:fig}
\end{figure}

Short computational experiments with the fine-scale model are used to construct
the coarse time-stepper -- the basic element for exchanging dynamical information
between coarse-scale model states and fine-scale states. 
A coarse time-stepper consists of three components: 
{\it lifting}, {\it fine-scale evolution} and {\it restriction} \cite{Gear:02} (Fig. \ref{cts:fig}).
The lifting transformation converts coarse-scale observables to consistent fine-scale states;
restriction is the reverse transformation, from fine-scale
sates to coarse-grained observables.  
Lifting, followed by restriction, should then give the identity on the coarse
observables (modulo roundoff error).
Different coarse time-steppers are generated via different lifting and 
restriction operators; one should test on-line that macroscopic 
computational results are insensitive to the specific 
details of the time-stepper implementation; a more extensive discussion can
be found in \cite{Kevrekidis:03}.

For many multiscale problems observed along a single effective spatial dimension, 
particle positions constitute the fine-scale model state,
while an obvious coarse-grained state is 
the local mean concentration of the particles \cite{Majda:99}.
For identical particles
this local mean concentration is observed in terms of the histogram of the single particle
position  probability density function (PDF).
However, this PDF histogram depends on the bin size used to estimate it; the PDF at any
given point becomes zero if the bin size is too small so that no particles exist 
within the bin containing this point. 
To overcome this difficulty, the cumulative distribution function (CDF) 
is naturally used as an alternative coarse-scale observable.  
The CDF has in principle infinite support, but its inverse, 
the Inverse CDF (ICDF), is supported in $[0,1]$; it can be readily represented 
by its projection on shifted Legendre polynomials \cite{Stegun:70}.

In multidimensional random particle systems, the inverse form of the multidimensional 
CDF is not readily available (strictly speaking, it is not even defined). 
The marginal and conditional ICDFs constitute candidate
coarse-scale observables for multidimensional systems; this is
because the multidimensional CDF can be represented using the 
marginal CDF in one direction and a collection of
conditional CDFs in the remaining directions. 
For instance, for systems in two spatial dimensions,
\begin{equation}
    F_{XY}(x,y)=\int_{-\infty}^y F_{X|Y}(x|y) {{dF_Y} \over {dy}}(y) dy,
\label{eqn1:eqn}
\end{equation}
where $F_{XY}(x,y)$, $F_Y(y)$ and $F_{X|Y}(x|y)$ are the CDF, 
differentiable marginal CDF and conditional CDFs, respectively. 
The conditional CDF, $F_{X|Y}(x|y)$, is defined by

\begin{eqnarray}
   F_{X|Y}(x|y) &=& \lim_{\Delta y \rightarrow 0} P (X \le x | y < Y \le y+\Delta y)   \nonumber \\
 &=& \lim_{\Delta y \rightarrow 0} {{F_{XY}(x,y+\Delta y)-F_{XY}(x,y)} \over {F_Y(y+\Delta y) - F_Y(y)}}.
\label{eqndefccdf:eqn}
\end{eqnarray}
Assuming smoothness, a finite number of conditional CDFs 
can be used to recover the particle distribution (e.g. through interpolation).
In the following, we illustrate through a two-dimensional system (without loss of
generality) our implementation of this procedure for multidimensional systems.

{\bf 2.1. Lifting.} 
Starting with the inverse CDFs (ICDF) for the marginal and conditional distributions, 
the lifting procedure
involves obtaining compatible realizations of the fine-scale states. 
Let the marginal ICDF in direction $y$, $IF_Y(\cdot): [0,1] \mapsto \mathbb R$, 
and conditional ICDFs in the other direction $x$, $IF_{X|Y}(\cdot,y): [0,1] \mapsto \mathbb R$, 
$y \in \mathbb R$, be defined by
\begin{eqnarray}
   IF_Y(f)=\arg_{y \in \mathbb R} \{ F_Y(y)=f \}, \quad f \in [0,1],  \nonumber
                                                                \\
   IF_{X|Y}(f,y)=\arg_{x \in \mathbb R} \{ F_{X|Y}(x|y)=f \}, \quad f \in [0,1]. 
\end{eqnarray}
First, the $y$-direction position of the $i^{th}$ particle is directly 
taken from the marginal ICDF as $y^s_i=IF_Y((i-0.5)/N),i=1,2,\cdots,N$, \cite{Gear:01} 
where $N$ is the number of particles and the superscript $s$ indicates that {$y^s_i$} 
is a sorted, monotonically ascending sequence. 
Then, corresponding to each $y^s_i$ generated in this manner, 
the $x$-direction position of the $i^{th}$ particle is determined 
as $x_i=IF_{X|Y}(U_i,y^s_i)$, where $U_i$ are i.i.d. real 
random variables with uniform distribution over $[0,1]$.
We only have a finite number of conditional CDFs, and we will assume that they are 
smooth in $y$; therefore, for each particular working {$y^s_i$} we employ the
conditional CDFs available in its neighborhood (e.g. the closest one, or possibly
an interpolation of the closest ones).
Only a few conditional ICDFs are needed if the CDF is sufficiently smooth. 
For example, if $M$ $(M \ll N)$ conditional ICDFs are needed, 
then these ICDFs, $IF_{X|Y}(f,y^c_k),k=1,2,\cdots,M,f \in [0,1]$, 
can be chosen such that $y^c_k=y^s_{(k-1)\cdot \mbox{int}(N/M)+\mbox{int}(N/2M)}$, 
where $\mbox{int}(\xi),\xi \in \mathbb R$ is the maximum integer not greater than $\xi$.    

{\bf 2.2. Fine-Scale Evolution.}
We work with multiscale dynamical systems for 
which we have {\it fine-scale} evolution models available. 
For noninteracting particles, letting $\bm X(t) = (x(t),y(t))^T \in{\mathbb R}^2$ 
denote the fine-scale model state, consisting of particle 
positions at time $t$, the discrete dynamics for $\bm X$ are given by

\begin{equation}
\bm X_{k+1}=\Psi(\bm X_{k},\bm \eta_k;\Delta t;\lambda), \quad k=0,1,\ldots
\label{disdyn:eqn}
\end{equation} 
where $\bm \eta_k\in {\mathbb R}^2$ denotes an 
external stochastic driving force and $\lambda$ 
is a set of (constant) parameters. 
In the class of problems envisioned here, the models $\Psi(\cdot)$ 
involve microscale simulators of Brownian motion, 
kinetic Monte Carlo simulation or molecular dynamics. 
    
{\bf 2.3. Restriction.}
Let the position of the $i^{th}$ particle be denoted as $(x_i,y_i),i=1,2,\cdots,N$ 
and the sorted particle positions in two directions be denoted as $\{ x^s_j \}$ and $\{y^s_m\}$, 
respectively. 
Then a two-dimensional mesh can be formed with each grid point 
having a coordinate $(x^s_j,y^s_m),j,m=1,2,\cdots,N$. 
For each point $(x^s_j,y^s_m)$, the number, $N_f$, of particles 
whose $x$ and $y$ direction positions satisfy $x_i \le x^s_j$ and 
$y_i \le y^s_m$, respectively, is counted and the CDF at 
this grid point evaluated as $F_{XY}(x^s_j,y^s_m)={{N_f-0.5} \over N}$. 
The CDF of particle positions can thus be obtained.
This is only one of several possible restriction methods;
other restriction approaches can be found in \cite{LiA:98}.

Assuming the CDF $F_{XY}(x,y)$ to be differentiable, 
Eqn. (\ref{eqn1:eqn}) or (\ref{eqndefccdf:eqn}) lead to the following fomula for the conditional CDF, 

\begin{equation}
 F_{X|Y}(x|y)={ {{{\partial F_{XY}} \over {\partial y}}(x,y)} \over { {{d F_Y} \over {d y}}(y)} }.
\label{eqn2:eqn}
\end{equation}
Continuing the example at the end of Section 2.1, 
the evaluation of the conditional ICDFs, $IF_{X|Y}(f,y^c_k),k=1,2,\cdots,M,f \in [0,1]$, 
requires the availability of the conditional CDFs, $F_{X|Y}(x|y^c_k)$.  
These can be numerically approximated from equation (\ref{eqn2:eqn}) as,

\begin{equation}
 F_{X|Y}(x^s_j|y^c_k)= {{F_{XY}(x^s_j,y^s_{p_2})-F_{XY}(x^s_j,y^s_{p_1})}\over{F_{XY}(x^s_N,y^s_{p_2})-F_{XY}(x^s_N,y^s_{p_1})}},
\label{eqn3:eqn}
\end{equation}   
where $p_1$ and $p_2$ can be chosen as $(k-1)\cdot \mbox{int}(N/M)+1$ 
and $k\cdot \mbox{int}(N/M)$, respectively. 
Once the conditional CDF $F_{X|Y}(x|y^c_k)$ is available numerically, 
the conditional ICDF $IF_{X|Y}(f,y^c_k)$ can be numerically evaluated 
as in the case of one-dimensional observables.

\vspace{0.5cm}
\noindent{\bf 3. Coarse Projective Integration (CPI) for Multidimensional Random Particle Systems}
\vspace{0.2cm}

\begin{figure*}
\centerline{\epsfig{figure=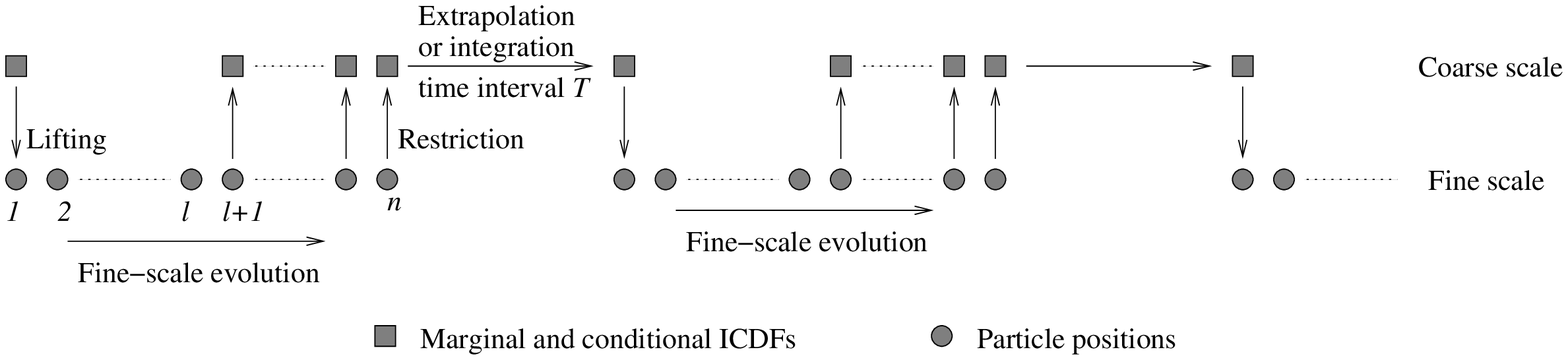,height=40mm,width=150mm,angle=0}}
\caption{The evolution of the coarse and fine-scale model states.}
\label{coarsefine:fig}
\end{figure*}

Coarse projective integrators (CPI) typically consist of four steps (Fig. \ref{coarsefine:fig}). 
At first, coarse observables are identified to which the fine-scale model states are restricted.
The coarse observables used in our context consist of the marginal and (finitely many) conditional 
ICDFs of the microscale particle positions, as described in the previous section; in particular, we use
a finite number of expansion coefficients (in some appropriate basis) 
of this marginal and these conditional ICDFs.
The particle positions can be generated through the {\it lifting} procedure 
described in Section 2.1 once an initial condition for these ICDFs 
have been specified. 
We represent the {\it lifting} operator by $\mu$, a mapping from the coarse 
observables $IF$ (ICDFs) to the microscopic descriptors ${\bm X}$ (particle positions). 
The second and third steps are the fine-scale evolution and 
restriction mentioned in sections 2.2 and 2.3, respectively.  
The restriction operator $\cal M$ is a mapping from the microscopic descriptors 
$\bm X$ to the coarse
observables $IF$, i.e., $IF={\cal M} {\bm X}$. 
Evidently, the operators $\cal M$ and $\cal \mu$ satisfy the property 
${\cal M}{\cal \mu}=I$ (modulo roundoff error).  
Along with restriction
comes the estimation of the coarse-scale time derivatives of the observables (the
marginal and ICDF coefficients).

The last step is the projection step in time -- the temporal evolution of our representation of the coarse-scale observables.
This step is templated on continuum numerical integration techniques - for coarse forward
Euler it is simple linear extrapolation of the coarse observables in time, although
more sophisticated and even implicit techniques can be (and have been) used
\cite{Gear:02,Rico-Martinez:04}.

Let the coarse-scale observables at time $t$ consist of 
$M+1$ ICDFs, $IF_{i,t},i=1,\cdots,M+1$, of which the first 
one is the marginal ICDF, $IF_Y(f)$, and the remaining ones are 
the conditional CDFs, $IF_{X|Y}(f,y^c_k)$, at $y$-direction positions $y^c_k,k=1,2,\cdots,M$. 
Let a basis of the coarse-scale subspace, 
which can be specified globally \cite{Gear:01,Setayeshgar:03} 
or locally \cite{Setayeshgar:03}, be denoted by $\{ \bm{\theta}_q \}$.
Then the projections, $\beta_{i,q,t}$, of ICDFs onto the basis can be computed by,

\begin{equation}
\beta_{i,q,t}=\left(IF_{i,t},\bm{\theta}_q \right), \quad q=0,1,\ldots,P \quad t=1,2,\cdots,n,
\label{eqn4:eqn}
\end{equation}
where the inner product is in $L_2$. 

The projective integration step over a coarse-scale time step $T$ 
can be formally written as

\begin{equation}
\beta_{i,q,n+T}=L(\beta_{i,q,l+1},\beta_{i,q,l+2},\cdots,\beta_{i,q,n}),
\label{eqn5:eqn}
\end{equation}
where $L(\cdot)$ is an operator based on (templated on) traditional continuum
numerical integration schemes, 
$\beta_{i,q,n+T}$ is the coefficient of the $q^{th}$ 
mode of the $i^{th}$ ICDF immediately after the temporal projection step 
and $\beta_{i,q,t},t=l+1,\cdots,n$ is the coefficient 
of the $q^{th}$ mode corresponding to the $t^{th}$ fine-scale time step prior to the projective step. 

Immediately after the projective sep, 
the new ICDFs, $IF_{i,n+T}$, based on the new coefficients $\beta_{i,q,n+T}$, are constructed as,

\begin{equation}
IF_{i,n+T}=\sum_{q=0}^P \beta_{i,q,n+T} \bm{\theta}_q \quad .
\label{eqn6:eqn}
\end{equation}
New fine-scale model states can then be ``lifted" from the ICDFs, $IF_{i,n+T}$.

Using the above steps, the procedure for applying CPI to multidimensional 
random particle systems can be summarized as follows:

\begin{enumerate} 
\item Generate fine-scale model state(s) consistent with the coarse-scale 
description given by the particle ICDFs.
The marginal and conditional ICDFs can be obtained 
through equation (\ref{eqn2:eqn}) if an analytical two-dimensional CDF is given as the initial condition. 
\item Let the fine-scale model state evolve according to the discrete dynamical model (\ref{disdyn:eqn}). 
\item Generate ICDFs at some successive fine-scale time steps. 
\item Project the ICDFs onto an appropriate basis (equation (\ref{eqn4:eqn})), 
and estimate the temporal derivatives of coefficients of the dominant modes.
\item Extrapolate (project forward in time) coefficients 
of the dominant modes over a large coarse-scale time interval 
$T$ (equation (\ref{eqn5:eqn})), reconstruct the ICDFs (equation (\ref{eqn6:eqn})) 
and go back to step 1.
\end{enumerate}
Usually one wants to report the multidimensional CDFs, and
they can be generated numerically along with the ICDFs 
using the approach in Section 2.3 
(although only the ICDFs are taken as coarse-scale states involved in the CPI method). 
The mesh size for numerically computing the CDF 
can be set larger than $\max_{1\le i \le N-1} (x^s_{i+1}-x^s_i,y^s_{i+1}-y^s_i)$,
in an attempt to alleviate fluctuation-related problems in the estimation step;
variance reduction schemes (multiple realizations of the simulation, or
more sophisticated approaches) may become necessary for this purpose.

\vspace{0.5cm}
\noindent{\bf 4. Coarse Dynamic Renormalization (CDR) for Multidimensional 
Random Particle Systems}
\vspace{0.2cm}

For multiscale systems of practical interest, if the
PDEs in the macroscopic level are scale invariant, 
they may possibly possess self-similar solutions \cite{Barenblatt:96}.
The analogy with traveling wave solutions for problems with translational
invariance is instructive: approaching a traveling wave in a co-traveling frame
appears like the approach to a stationary state.
Similarly, approaching a self-similar solution in a dynamically renormalized
(co-exploding or co-collapsing) frame, appears like the approach to a stationary state.
Dynamic renormalization procedures have been used to investigate  
self-similar systems \cite{Mclaughlin:86,LemesurierA:88,LemesurierB:88}; 
recently a template-based
approach for studying the dynamics of problems with translational 
symmetry \cite{Rowley:00}
has been extended to study the dynamics of problems with scale invariance
\cite{Aronson:01,Siettos:03,Rowley:03,Chen:03}.
When macroscopic scale-invariant PDEs are explicitly available, template
conditions can be used to derive dynamical equations (termed ``MN-dynamics")
for the rescaled self-similar solutions and similarity exponents \cite{Aronson:01}. 
The idea of employing template conditions can also be used to obtain renormalized 
self-similar macroscale solutions and similarity exponents for multiscale systems whose coarse-level
PDEs are not explicitly known \cite{Chen:03}. 
The number of template conditions depends on how many 
rescaling variables are needed to renormalize the physical solutions.
   
Consider a PDE in the form of

\begin{equation}
  {{\partial F} \over {\partial t}}=D_{xy}(F),
\label{eqnoriginal:eqn}
\end{equation}
where $F(x,y,t)$ is a CDF of particle positions which do not 
collectively translate in the space domain
(the case of joint scale and translational invariance can also be simply treated,
see \cite{Rowley:03}).
  
The differential operator $D_{xy}$ is such that there exist constants $p$ and $a$ such that

\begin{equation}
   D_{xy}(f({x \over A},{y \over {A^p}}))=A^a D_{uv}(f(u,v)),\quad u= {x \over A}, \quad v={y \over {A^p}},
\label{eqnDxy:eqn}
\end{equation}
for any real function $f$, real value $A>0$ and coordinate $(x,y)$ (there is
no amplitude rescaling since this is a CDF).
If a  self-similar solution $F(x,y,t)$ exists, it can be written as 

\begin{equation}
   F(x,y,t)=U({x \over (cs)^\alpha}, {y \over (cs)^{\alpha p}};c),
\label{eqn7:eqn}
\end{equation}
where $c$ is a constant parametrizing the family of 
self-similar shapes, $s=t-t_0$, ($t_0$ is the blowup time for 
problems with finite time singularities) and $t>t_0$.
Substituting (\ref{eqn7:eqn}) into (\ref{eqnoriginal:eqn}), we have  

\begin{equation}
  \alpha a = -1,
\label{eqnalphaa:eqn}
\end{equation}
and $U$ satisfies the PDE,

\begin{equation}
  - \alpha u U_u- \alpha p v U_v = c^{-1} D_{uv}(U),
\label{eqnselfsimilar:eqn}
\end{equation}
where $u=x/(cs)^\alpha$, $v=y/(cs)^{\alpha p}$. For the operator $D_{xy}$ that satisfies Eqn. (\ref{eqnDxy:eqn}), the constant $a$ is determined by $D_{xy}$ itself. Hence the similarity exponent $\alpha$ can be calculated by Eqn.(\ref{eqnalphaa:eqn}).

If the macroscopic equation (\ref{eqnoriginal:eqn}) is not explicitly known, 
one cannot analytically obtain
the exponents $p$ and $a$; numerical computations are needed to determine these
constants -- and thus to test the scale invariance of the operator -- before
locating the self-similar solutions themselves.
For an operator $D_{xy}$ that satisfies the (unknown) equation (\ref{eqnDxy:eqn}), 
the constants $p$ and $a$ can be obtained using a black box simulator of the
equation as follows:
Since the unknown Eqn. (\ref{eqnDxy:eqn}) is valid for {\it any} coordinate $(x,y)$ and real function $f$, let $f$ be a test function
(we choose it here for convenience to be exponential in space), $(x_1,y_1)=(u_1 A,v_1 A^p)$ and $(x_2,y_2)=(u_2 A,v_2 A^p)$, 
where $A$ is arbitrarily chosen as a positive real value.
This would imply that (choosing two points in space) the following two relations hold:

\begin{eqnarray}
&&D_{xy}(f({x \over A},{y \over {A^p}}))(x_1,y_1)= A^a D_{uv}(f(u,v))(u_1,v_1),  \nonumber
                                                                              \\
&&D_{xy}(f({x \over A},{y \over {A^p}}))(x_2,y_2)= A^a D_{uv}(f(u,v))(u_2,v_2),  \nonumber
\end{eqnarray}
where $(u_1,v_1)$ and $(u_2,v_2)$ are two distinct coordinates.
Comparing the above two equations, we have

\begin{equation}
   {{D_{xy}(f({x \over A},{y \over {A^p}}))(x_1,y_1)} \over {D_{xy}(f({x \over A},{y \over {A^p}}))(x_2,y_2)}}= {{D_{uv}(f(u,v))(u_1,v_1)} \over {D_{uv}(f(u,v))(u_2,v_2)}}.
\label{eqnconp:eqn}
\end{equation}
Therefore, the constant $p$ is the solution to Eqn. (\ref{eqnconp:eqn}). 
We rewrite Eqn.(\ref{eqnconp:eqn}) as 
\begin{eqnarray}
& &   D_{xy}(f({x \over A},{y \over {A^p}}))(x_2,y_2) \\ \nonumber 
&- &{{D_{xy}(f({x \over A},{y \over {A^p}}))(x_1,y_1)} \over {{D_{uv}(f(u,v))(u_1,v_1)}}} D_{uv}(f(u,v))(u_2,v_2)=0. 
\label{eqnpsolve:eqn}
\end{eqnarray}
The constant $p$ can then be solved for, employing an in principle
arbitrary test function $f$, and 
using Newton's method. 
Since we assumed that the operator $D_{xy}$ is not explicitly available, 
$D_{xy}(f(x,y))$ can be estimated
by running the micro-simulator for short time ``bursts"
and numerically obtaining the derivative ${\partial f} \over {\partial t}$. 
The constant $a$ is calculated by 
\begin{equation}
    a=\log_{A}{{D_{xy}(f({x \over A},{y \over {A^p}}))(x_1,y_1)} \over {{D_{uv}(f(u,v))(u_1,v_1)}}},
\end{equation} 
once $p$ is obtained.
Clearly, other test functions and conditions evaluated at other points 
 can be used;
care must be taken also to ensure the finiteness of the estimated quantities.

Given $p$, to determine the self-similar shape of the solution, 
we consider the general scaling
\begin{equation}
  F(x,y,t)=\omega ( {x \over A(t)}, {y \over {A(t)^p}}, t),
\label{eqnrescaling:eqn}
\end{equation}
where $A(t)$ is an unknown function. 
The PDE becomes
\begin{equation}
  \omega_t-{ A_t \over A} u \omega_u - {{p A_t} \over A} v \omega_v = A^a D_{uv}(\omega).
\label{eqnomega:eqn}
\end{equation}
Evidently, $U$ and $\omega$ are both renormalized CDFs.

By comparing equations (\ref{eqnselfsimilar:eqn}) and (\ref{eqnomega:eqn}), we have 
\begin{equation}
   \alpha = {{\lim_{t \rightarrow \infty} { A_t \over A}} \over {\lim_{t \rightarrow \infty} c A^a}}.
\label{alphaA:eqn}
\end{equation}

Also by comparing equations (\ref{eqn7:eqn}) and (\ref{eqnrescaling:eqn}), we have 
\begin{equation}
   {{\lim_{t \rightarrow \infty} A} \over {\lim_{t \rightarrow \infty} (cs)^{\alpha}}} =1.
\label{At:eqn}
\end{equation}

The above equations (\ref{alphaA:eqn}) and (\ref{At:eqn}) together with (\ref{eqnalphaa:eqn}) lead to 
\begin{equation}
   \alpha = {\lim_{t \rightarrow \infty} { {(t-t_0) A_t} \over A}}.
\label{alphalim:eqn}
\end{equation}

Therefore, the value for $\alpha$ can be calculated once values of 
$A(t)$ are obtained in the long-time limit (i.e.,
after $\omega$ reaches the steady state). 
Indeed, let $t_1$ and $t_2$ be distinct times 
after $\omega$ reaches the steady state, then by  
(\ref{alphalim:eqn}),
\begin{equation}
   \alpha= { {t_2-t_1} \over { { A(t_2) \over A_t(t_2)} -{ A(t_1) \over A_t(t_1)}} }.
\label{alphacal:eqn}
\end{equation}

A single template condition 
is required to solve for both $\omega(u,v,t)$ and  $A(t)$
at every time step. 
In this paper, the template is chosen to be
\begin{equation}
  \omega(e,\infty,t)=m, \quad e<0, \quad  0<m<0.5,
\end{equation}
where $e$ and $m$ are both constants. 
The template condition has the following physical meaning: 
the rescaled marginal CDF $\omega_U$ always has the 
same value $m$ at the u-direction coordinate $e$ for all time $t$.
Applying this template to Eqn.(\ref{eqnomega:eqn}) and 
assuming ${{\partial \omega} \over {\partial v}}(e,v,t)$ decays 
exponentially as $v \rightarrow \infty$, we have
\begin{equation}
  A_t e {{\partial \omega} \over {\partial u}} (e,\infty,t) + A^{a+1} D_{uv}(\omega)(e,\infty,t) = 0.
\label{eqncouple:eqn}
\end{equation}  

\begin{figure*}
\centerline{\epsfig{figure=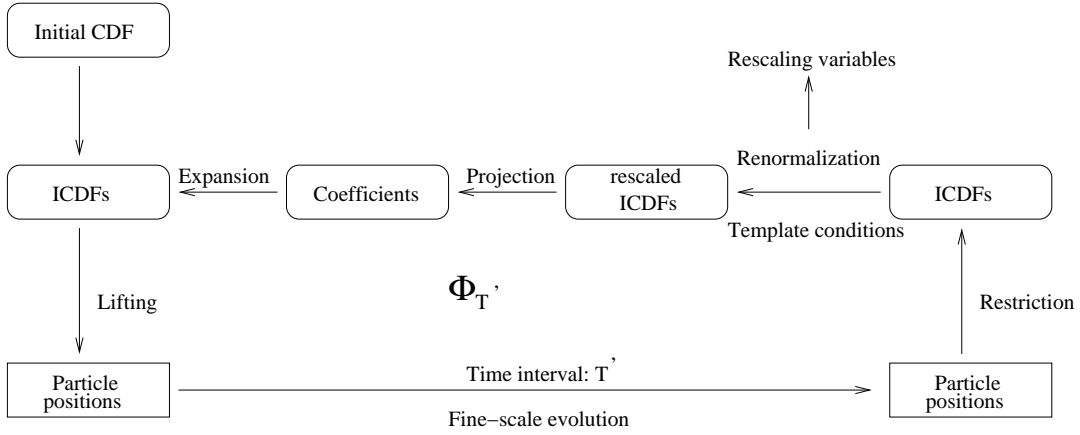,width=0.8\textwidth}}
\caption{Schematic illustration of coarse dynamic renormalization.}
\label{crm:fig}
\end{figure*}

Equations (\ref{eqnomega:eqn}) and (\ref{eqncouple:eqn}) can be coupled to 
solve for the rescaled CDF $\omega$ and rescaling variable $A$ if the 
operator $D_{xy}$ is explicitly known.
As the time $t \rightarrow \infty$, $\omega$ may approach a steady state, 
which is then a stable self-similar shape for the 
solutions to Eqn.(\ref{eqnoriginal:eqn}). 
In general cases, the macroscale equation for the 
CDF of the particle positions may not be explicitly available.
However, the template-based approach can still be used to 
renormalize the CDF evolved via microscale models and 
rescaling variables are obtained during the course of renormalization. 
In these cases, to express the CDFs, the marginal and conditional 
ICDFs are used again as macroscopic observables,
and their (discretized) projections over an orthonormal basis are again 
used to numerically characterize them through a finite number of
coefficients. 
Based on the coarse time-stepper, the procedure for the 
coarse renormalization is schematically depicted in Fig. \ref{crm:fig} 
and consists of the following steps, 
\begin{enumerate}
\item Generate the marginal and conditional ICDFs $IF_{i,t},i=1,\cdots,M+1$ 
according to the initial CDF using equation (\ref{eqn2:eqn}) 
or (\ref{eqn3:eqn}) or according to the 
coefficients, $\beta_{i,q,t}$, of dominant modes 
of the ICDFs using equation (\ref{eqn6:eqn}).
\item Generate particle positions using the ICDFs 
using the {\it lifting} procedure in the coarse time-stepper.
\item Evolve particle positions over a time 
interval $T'$ using a fine-scale model (\ref{disdyn:eqn}).
\item Obtain ICDFs from particle positions using 
the {\it restriction} procedure in the coarse time-stepper.
\item Rescale the marginal ICDF according to the 
template condition and obtain the rescaling variable $A$. 
We then rescale the conditional ICDFs by a factor of $A^p$. 
This step can be justified by Eqn.(\ref{eqnrescaling:eqn}). 
Indeed, obtaining the rescaled solution 
$\omega_{k+1}$ from $\omega_k$ via the dynamics 
(\ref{eqnomega:eqn}) and (\ref{eqncouple:eqn}) is 
equivalent to starting from the initial condition $\omega_k$, 
running the original dynamics for a while to get $F_{k+1}$, 
factoring out the rescaling variable $A$, and 
rescaling $F_{k+1}$ in scales of $A$ and $A^p$ 
respectively in $x$ and $y$ directions. 
\item Project the rescaled ICDFs onto 
the orthonormal basis and obtain the coefficients of 
leading modes using equation (\ref{eqn4:eqn}). Go back to step 2.
\end{enumerate} 

The above procedure can be viewed as an iterative 
algorithm to solve the fixed point of a nonlinear operator $\Phi_{T'}$,
written as,
\begin{equation}
{\bm \beta}=\Phi_{T'} ({\bm \beta}) .
\label{eqn8:eqn}
\end{equation}
This fixed point can be written in component form as,
$\beta^r_{i,q},i=1,\cdots,M+1,q=0,\cdots,P$, \cite{Chen:03} or
\begin{eqnarray}
{\bm \beta} &=&(\beta^r_{1,0},\beta^r_{1,1},\cdots,\beta^r_{1,P},\beta^r_{2,0},\beta^r_{2,1},\cdots,\beta^r_{2,P}, \nonumber \\
&\cdots&,\beta^r_{M+1,0},\beta^r_{M+1,1},\cdots,\beta^r_{M+1,P})^T,  \nonumber 
\end{eqnarray}
where the superscript $r$ refers 
to the fact that this is the renormalized self-similar shape. 
These coefficients correspond to the renormalized self-similar 
ICDFs and CDF of the multidimensional particle system. 
Equation (\ref{eqn8:eqn}) may be solved using any numerical 
algorithm such as direct iteration or matrix-free (Krylov-subspace based)
implementations of
Newton's method \cite{Kelley:95}.

\vspace{0.5cm}
\noindent{\bf 5. Self-Similar and Asymptotically Self-Similar Dynamics of
Brownian Particles in a Couette Flow}
\vspace{0.2cm}

We will use CPI and CDR algorithms to study 
two-dimensional Brownian models of particle dispersion
in a Couette flow \cite{Panton:96}. 
In this section, a particle system with self-similar dynamic
evolution and a system with asymptotically
self-similar evolution are explored, respectively. 

{\bf 5.1. A Particle System with Self-Similar Dynamics.}
Let $X(t)$ and $Y(t)$ represent particle positions 
in $x$ and $y$ directions respectively at time $t$ 
on an infinite two-dimensional spatial domain. 
The particle positions in the two directions evolve
in this model governed by the following dynamics:
\begin{equation}
  dX(t)=DdW_X (t), \quad dY(t)=Xdt,
\label{eqn9:eqn}
\end{equation}
where $W_X (t)$ is a Wiener processes \cite{Gihman:72} 
and $D$ is the diffusion coefficient. 
The discretized dynamics of (\ref{eqn9:eqn}) is given by \cite{Milshtein:74}
\begin{equation}
   \quad X_{k+1}=X_k+D\eta_{X,k}\sqrt{\Delta t}, \quad Y_{k+1}=Y_k+X_k \Delta t,
\label{eqn10:eqn}
\end{equation}
where $\eta_{X,k}$ are i.i.d. standard Gaussian random variables. 

The dynamics (\ref{eqn9:eqn}) represent the motion of 
particles which only diffuse in the $x$-direction in a Couette flow.
It can be shown that the coarse-scale dynamics for the PDF, 
$P_{XY}(x,y,t)$, of a particle position, 
corresponding to the fine-scale dynamics (\ref{eqn9:eqn}), 
is governed by the following equation \cite{Majda:99}

\begin{equation}
  {{\partial P_{XY}} \over {\partial t}} + x {{\partial P_{XY}} \over {\partial y}} = {{D^2}\over 2} {{\partial^2 P_{XY}} \over {\partial x^2}},
\label{eqn11:eqn}
\end{equation} 
where $P_{XY}$ is assumed to be 2nd-order differentiable. 
Hence the dynamics for the CDF, $F_{XY}(x,y,t)$, 
associated with (\ref{eqn11:eqn}) is given by
\begin{equation}
  {{\partial F_{XY}} \over {\partial t}} + x {{\partial F_{XY}} \over {\partial y}} - \int_{-\infty}^x {{\partial F_{XY}} \over {\partial y}} dx = {{D^2}\over 2} {{\partial^2 F_{XY}} \over {\partial x^2}} \ .
\label{eqn12:eqn}
\end{equation}

In the above equation, the operator $D_{xy}$ is written as
\begin{equation}
  D_{xy}= - x {{\partial } \over {\partial y}} + \int_{-\infty}^x {{\partial} \over {\partial y}} dx + {{D^2}\over 2} {{\partial^2 } \over {\partial x^2}} .
\end{equation}

This operator satisfies the scale invariance property (\ref{eqnDxy:eqn}) for constant values $p=3$ and $a=-2$. 
The analytical self-similar solution to Eqn. (\ref{eqn11:eqn}) is inspired by \cite{Okubo:69} (See Appendix A)
\begin{equation}
  P_{XY}(x,y,t)={\sqrt{3} \over {\pi {D^2} (t-t_0)^2 }} e^{- \left( {{6(y - 0.5x (t-t_0))^2} \over {D^2 (t-t_0)^3}} + {{x^2} \over {2D^2 (t-t_0)}} \right)}, 
\label{eqn102:eqn}
\end{equation} 
where $t_0$ is the blowup time 
(backward in time), which then gives the self-similar solution to (\ref{eqn12:eqn}),
\begin{eqnarray}
  F_{XY}(x,y,t) &=&{\sqrt{3} \over {\pi {D^2} (t-t_0)^2}}  \nonumber \\
\int_{-\infty}^x \int_{-\infty}^y &e&^{- ( {{6(y - 0.5x(t-t_0))^2} \over {D^2 (t-t_0)^3}} + {{x^2} \over {2D^2 (t-t_0)}} )} dy dx. 
\label{eqn103:eqn}
\end{eqnarray} 

Let $u'={x \over {(c(t-t_0))^{1/2}}}$ and $v'={y \over {(c(t-t_0))^{3/2}}}$, then 
\begin{eqnarray}
F_{XY}(x,y,t) &=& F_{UV}(u',v')={ {\sqrt{3}} \over {\pi D^2/c^2}} \nonumber \\
\int_{-\infty}^{u'} \int_{-\infty}^{v'} &e&^{- ( {{6(v - 0.5u/c)^2} \over {D^2/c^3}} + {{u^2} \over {2D^2/c}} )} dv du.
\label{eqn104:eqn}
\end{eqnarray}

Hence for the integro-differential equation (\ref{eqn12:eqn}), 
the similarity exponent in (\ref{eqn7:eqn}) is $\alpha=1/2$. 
For the CDF in (\ref{eqn104:eqn}), its standard deviations (std.'s) in two 
directions and correlation are $\sigma_X=D/c^{1/2}$, 
$\sigma_Y=D/(\sqrt{3}c^{3/2})$, and $\rho_{XY}=\sqrt{3}/2$, respectively. 

{\bf 5.2. A Particle System with Asymptotically Self-Similar Dynamics.}
We now consider Brownian particles in a Couette flow that diffuse
in both spatial directions.  
The microscopic-level evolution equation for the particle positions is given by
\begin{equation}
  dX(t)=DdW_X (t), \quad dY(t)=Xdt+DdW_Y(t),
\label{2ddifBrownian:eqn}
\end{equation}
where $W_X (t)$ and $W_Y(t)$ are independent Wiener processes. 
The discretized dynamics of (\ref{2ddifBrownian:eqn}) are given by \cite{Milshtein:74}
\begin{eqnarray}
   X_{k+1} &=&X_k+D\eta_{X,k}\sqrt{\Delta t}, \nonumber \\
   Y_{k+1} &=&Y_k+X_k\Delta t+D\eta_{Y,k}\sqrt{\Delta t},
\label{dis2ddifBrownian:eqn}
\end{eqnarray}
where $\eta_{X,k}$ and $\eta_{Y,k}$ are i.i.d. standard Gaussian random variables. 

The coarse-scale PDF of the particle positions is governed by  \cite{Majda:99}
\begin{equation}
  {{\partial P_{XY}} \over {\partial t}} + x {{\partial P_{XY}} \over {\partial y}} = {{D^2}\over 2} {{\partial^2 P_{XY}} \over {\partial x^2}}+{{D^2}\over 2} {{\partial^2 P_{XY}} \over {\partial y^2}}.
\label{pdf2ddifBrownian:eqn}
\end{equation} 
Accordingly, the dynamics for the CDF, $F_{XY}(x,y,t)$, are given by
\begin{eqnarray}
  {{\partial F_{XY}} \over {\partial t}} &+& x {{\partial F_{XY}} \over {\partial y}} - \int_{-\infty}^x {{\partial F_{XY}} \over {\partial y}} dx \nonumber \\
&=& {{D^2}\over 2} {{\partial^2 F_{XY}} \over {\partial x^2}} +{{D^2}\over 2} {{\partial^2 F_{XY}} \over {\partial y^2}} \ .
\label{cdf2ddifBrownian:eqn}
\end{eqnarray}

In Eqn. (\ref{cdf2ddifBrownian:eqn}), the operator $D_{xy}$ is written as
\begin{equation}
  D_{xy}= - x {{\partial } \over {\partial y}} + \int_{-\infty}^x {{\partial} \over {\partial y}} dx + {{D^2}\over 2} {{\partial^2 } \over {\partial x^2}}+{{D^2}\over 2} {{\partial^2 } \over {\partial y^2}} \ .
\end{equation}

The above operator does not exactly satisfy the scale invariance 
property (\ref{eqnDxy:eqn}) for any function $f$, 
positive real value $A$ and coordinate $(x,y)$. 
However, as {\it the scale} of the function $f$ or the value of $A$ 
become sufficiently large, 
there may exist $p$ and $a$ such that (\ref{eqnDxy:eqn}) is {\it approximately} 
satisfied. 
At those limits for $f$ and/or $A$, 
values of $p$ and $a$ that approximately 
satisfy (\ref{eqnDxy:eqn}) approach the values 
that exactly satisfy the scale invariance property in the true self-similar case.  

The asymptotically self-similar solution to the 
equation (\ref{pdf2ddifBrownian:eqn}) is given by \cite{Baptista:95}
\begin{eqnarray}
& & P_{XY}(x,y,t)={1 \over {2 \pi {D^2} (t-t_0) (1+(t-t_0)^2/12)^{1/2}}} \nonumber \\
& &e ^{- \left( {{(y - 0.5x(t-t_0))^2} \over {2D^2(t-t_0)(1+(t-t_0)^2/12)}} + {{x^2} \over {2D^2(t-t_0)}} \right)},
\label{pdf2ddifsolution:eqn}
\end{eqnarray} 
which then provides the asymptotically self-similar 
solution to (\ref{cdf2ddifBrownian:eqn}),
\begin{eqnarray}
 & & F_{XY}(x,y,t)={1 \over {2 \pi {D^2} (t-t_0) (1+(t-t_0)^2/12)^{1/2}}} \nonumber \\
& &\int_{-\infty}^x \int_{-\infty}^y e^{- ( {{(y - 0.5x(t-t_0))^2} \over {2D^2(t-t_0)(1+(t-t_0)^2/12)}} + {{x^2} \over {2D^2(t-t_0)}} )} dy dx. \nonumber \\
& &
\label{cdf2ddifsolution:eqn}
\end{eqnarray} 

Let $u'={x \over {(c(t-t_0))^{1/2}}}$ and $v'={y \over {(c(t-t_0))^{3/2}}}$, then 
\begin{eqnarray}
  & &F_{XY}(x,y,t) =F_{UV}(u',v',t) = \nonumber \\ 
& &{ {\sqrt{3}(t-t_0)} \over {\pi D^2 ((t-t_0)^2+12)^{1/2}/c^2}} \nonumber \\
& &\int_{-\infty}^{u'} \int_{-\infty}^{v'} e^{- ( {{6(v - 0.5u/c)^2 (t-t_0)^2} \over {D^2 ((t-t_0)^2+12)/c^3}} + {{u^2} \over {2D^2/c}} )} dv du.
\label{cdf2ddifrescaledsol:eqn}
\end{eqnarray}

In the long-time limit, $F_{UV}(u',v',t)$ has a steady-state form, 
which is the same as that given by (\ref{eqn104:eqn}) in the self-similar case. 

\vspace{0.5cm}
\noindent{\bf 6. Numerical Examples}
\vspace{0.2cm}

In what follows, direct particle simulations are implemented to 
accelerate the numerical evolution of the CDFs via CPI, and
to locate self-similar CDFs via CDR 
for the particle systems in Section 5. 
The approximate (asymptotic) scale invariance of the 
macroscale differential operator $D_{xy}$ for the 
asymptotically self-similar particle system is also examined. 
The fixed-point algorithm in Section 4 is 
utilized to solve for the long-time steady-state shape 
of the CDF for the asymptotically self-similar system. 
The diffusion coefficient $D$ and simulation time step $\Delta t$ 
in the fine-scale model are set to $5.0cm/s^{1/2}$ and $0.01s$, respectively. 
An ensemble of 2000 ($N=2000$) particles is used in the fine-scale 
simulations except where otherwise indicated. 

\bigskip

{\noindent \bf Simulation 1: Direct Simulation of the Self-Similar Particle System} 

\begin{figure}
\begin{minipage}{4cm}
\epsfig{figure=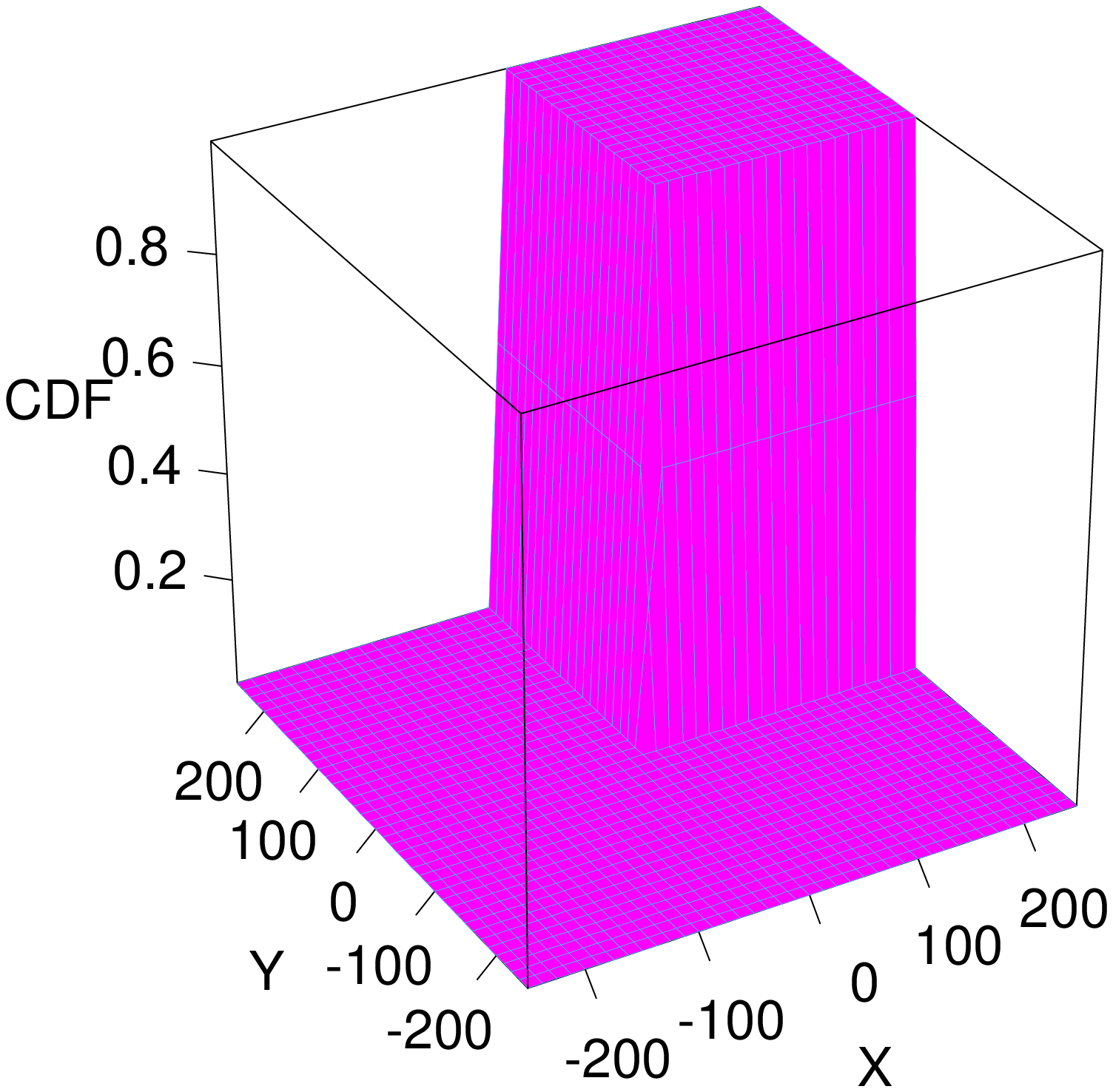,width=\textwidth}
\end{minipage}
\hfill
\begin{minipage}{4cm}
\epsfig{figure=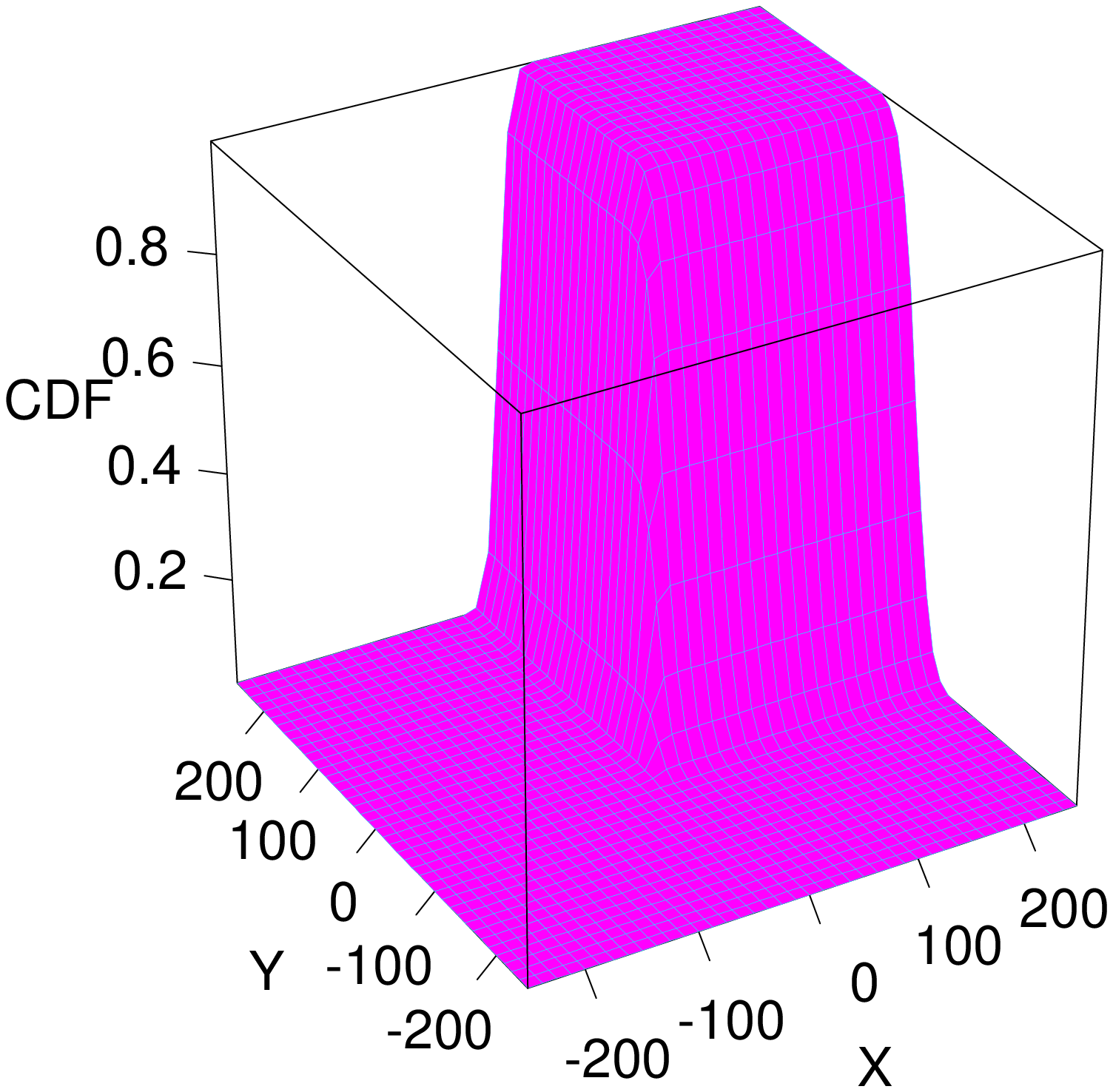,width=\textwidth}
\end{minipage}
\hfill
\begin{minipage}{4cm}
\epsfig{figure=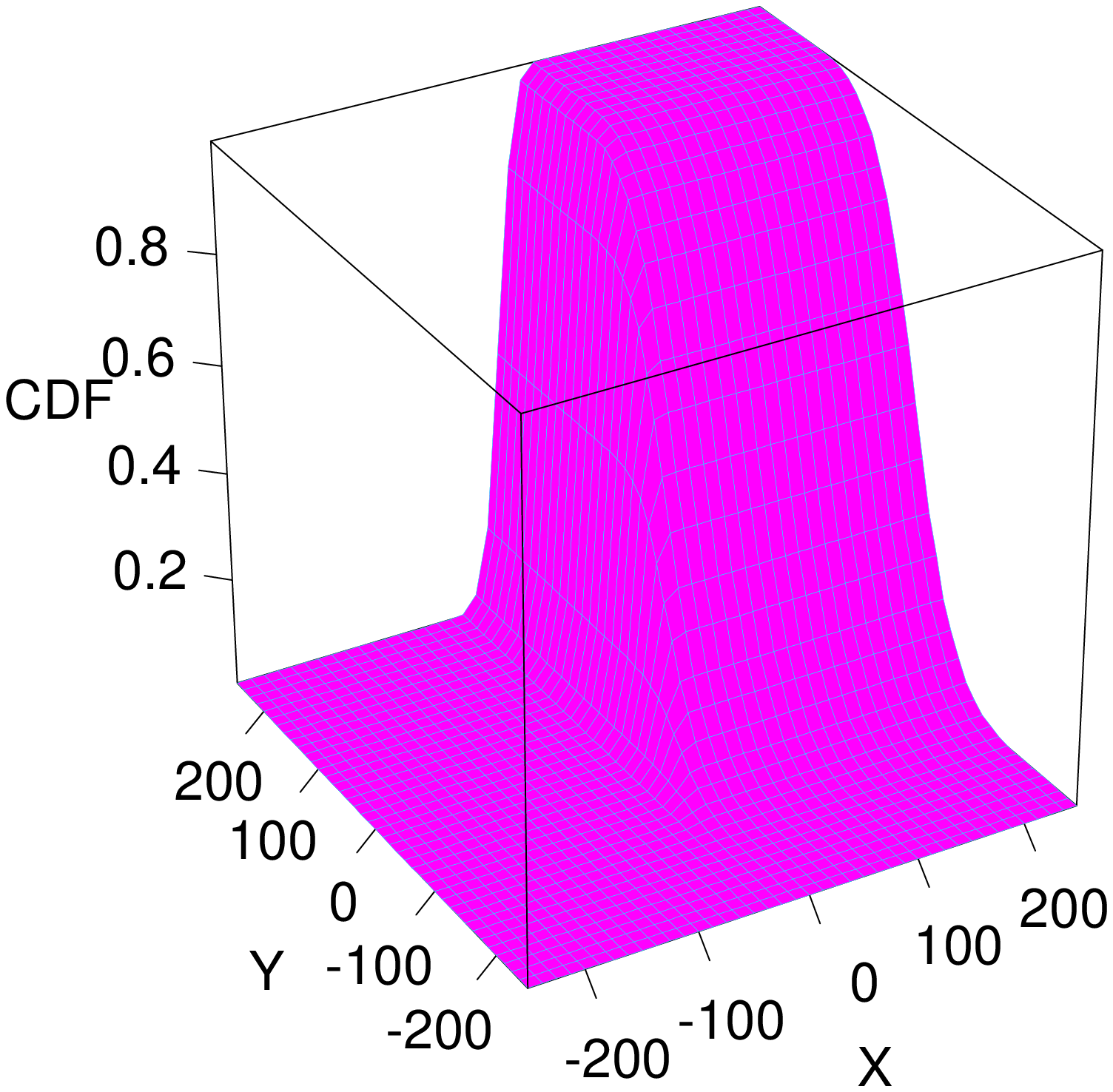,width=\textwidth}
\end{minipage}
\hfill
\begin{minipage}{4cm}
\epsfig{figure=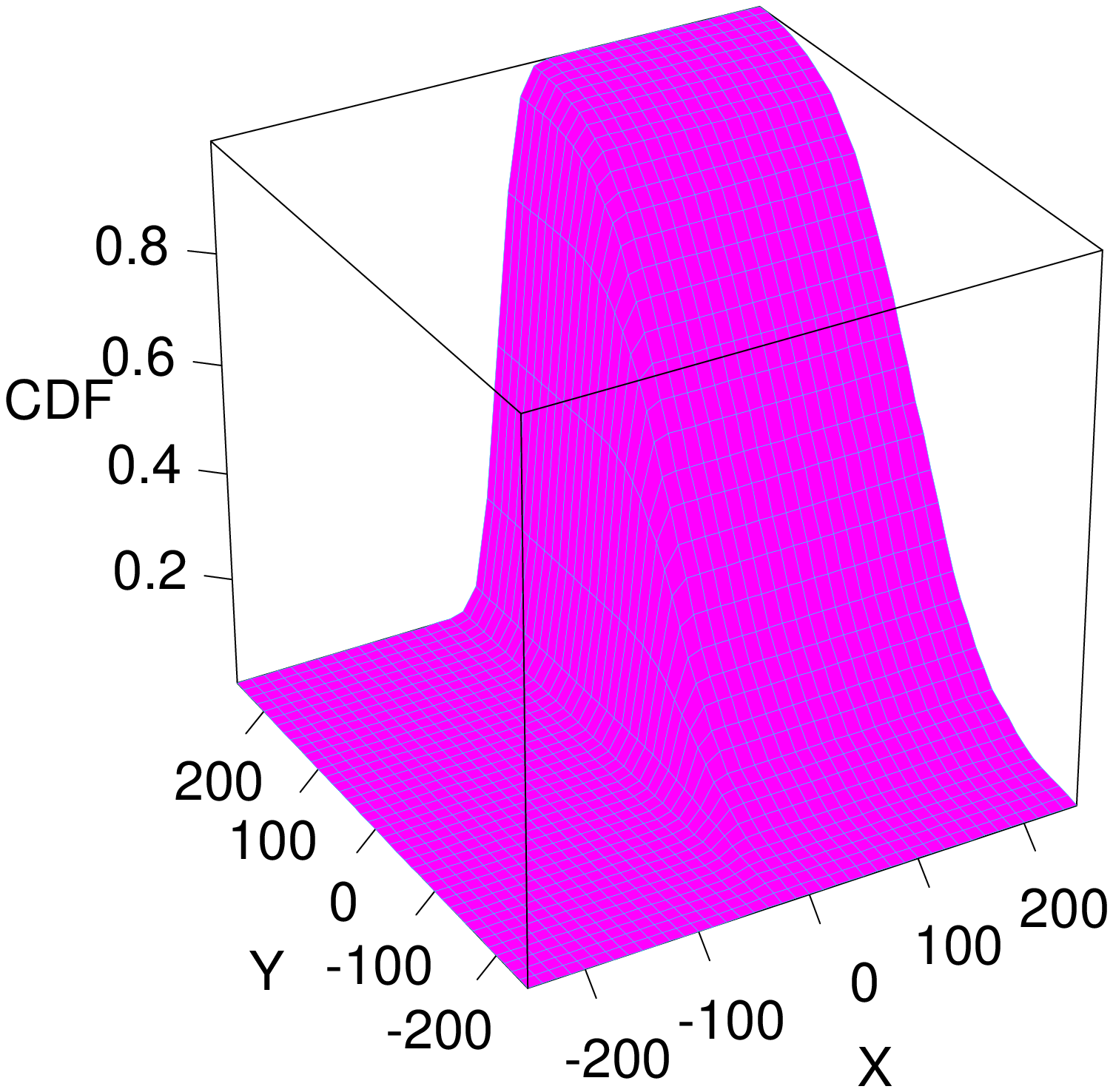,width=\textwidth}
\end{minipage}
\hfill
\caption{True CDFs at a sequence of time steps; top left: $t=0$, top right: $t=300 \Delta t$, bottom left: $t=600 \Delta t$, bottom right: $t=900 \Delta t$.} 
\label{fig4:fig}
\end{figure}

In this simulation, the initial fine-scale particle positions 
are chosen to follow a uniform distribution over the 
square domain $(-10cm,10cm)\times(-10cm,10cm)$.
%
%
An ensemble of 2000 particles whose distribution is consistent with
the coarse-scale initial conditions are evolved directly 
using (\ref{eqn10:eqn}) and used to construct true evolved coarse-scale CDFs. 
Particle positions are recorded at the time $300\Delta t$, 
$600 \Delta t$ and $900 \Delta t$, respectively. 
Two-dimensional CDFs are numerically computed using the 
procedure in Section 2.3 and plotted in Fig. \ref{fig4:fig}.
The number of grid points used to compute the CDFs is 1681. 

\bigskip
{\noindent\bf Simulation 2: Coarse Projective Integration of the 
Self-Similar Particle System}

\begin{figure}
\begin{minipage}{4cm}
\epsfig{figure=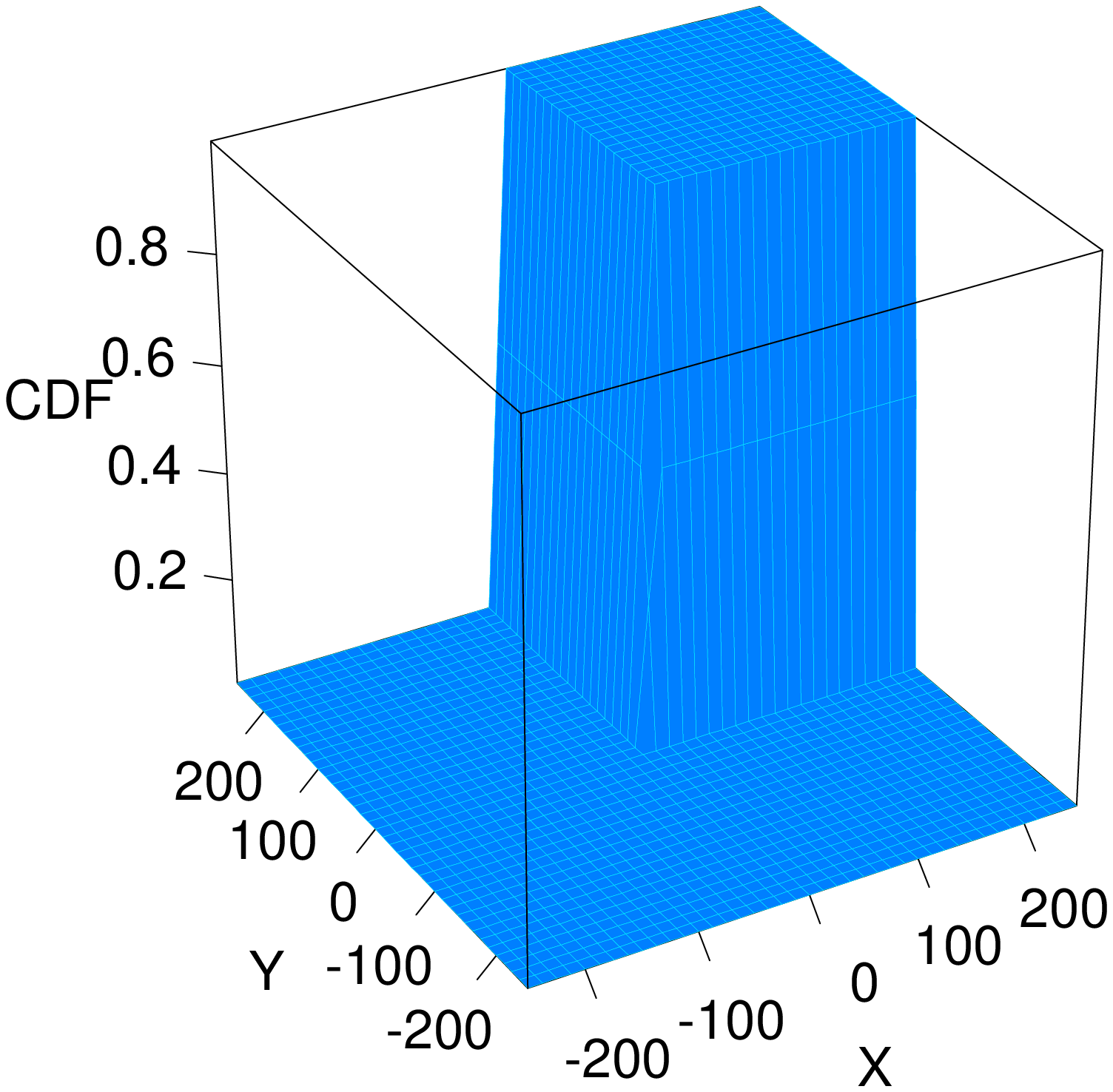,width=\textwidth}
\end{minipage}
\hfill
\begin{minipage}{4cm}
\epsfig{figure=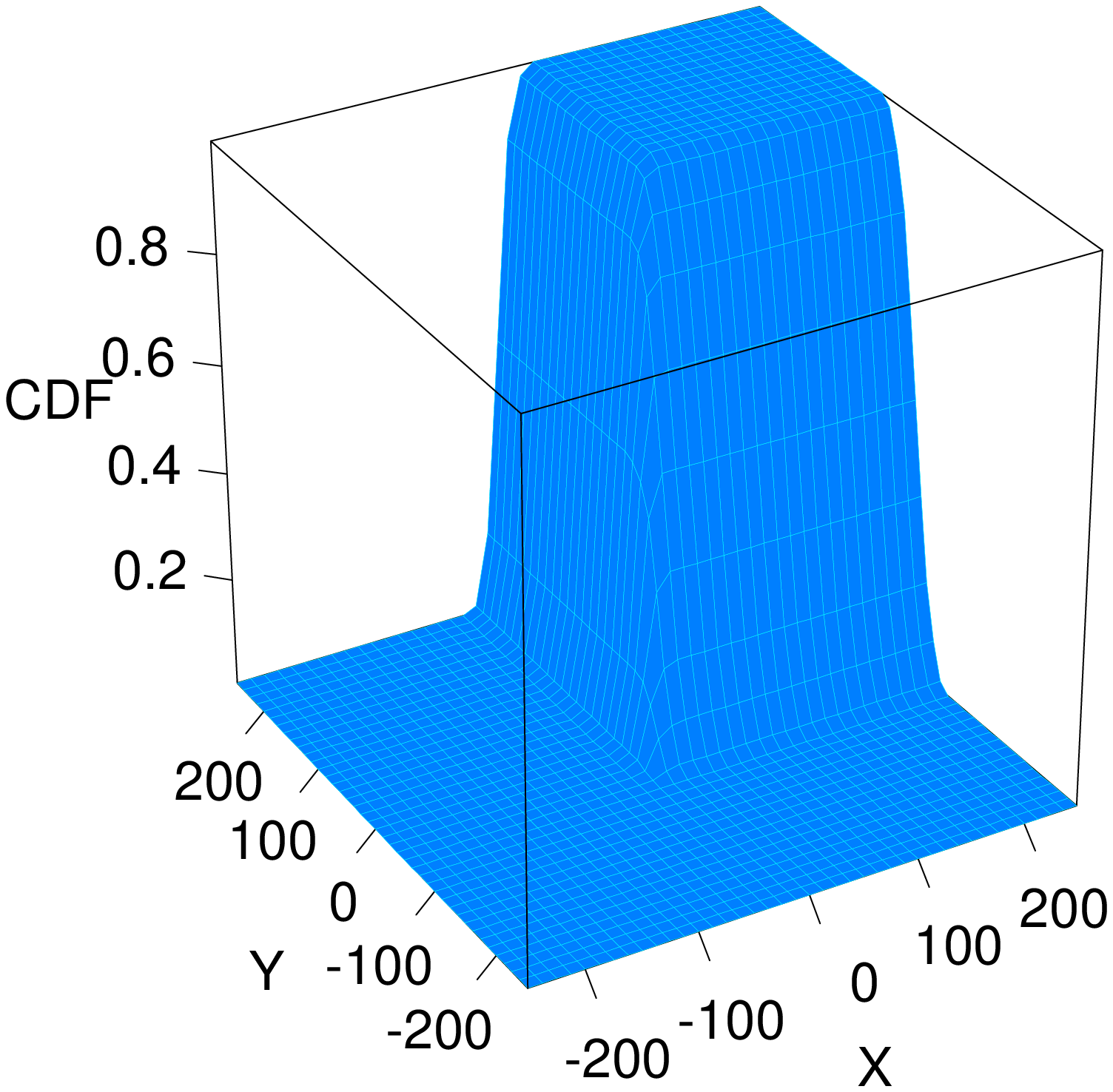,width=\textwidth}
\end{minipage}
\hfill
\begin{minipage}{4cm}
\epsfig{figure=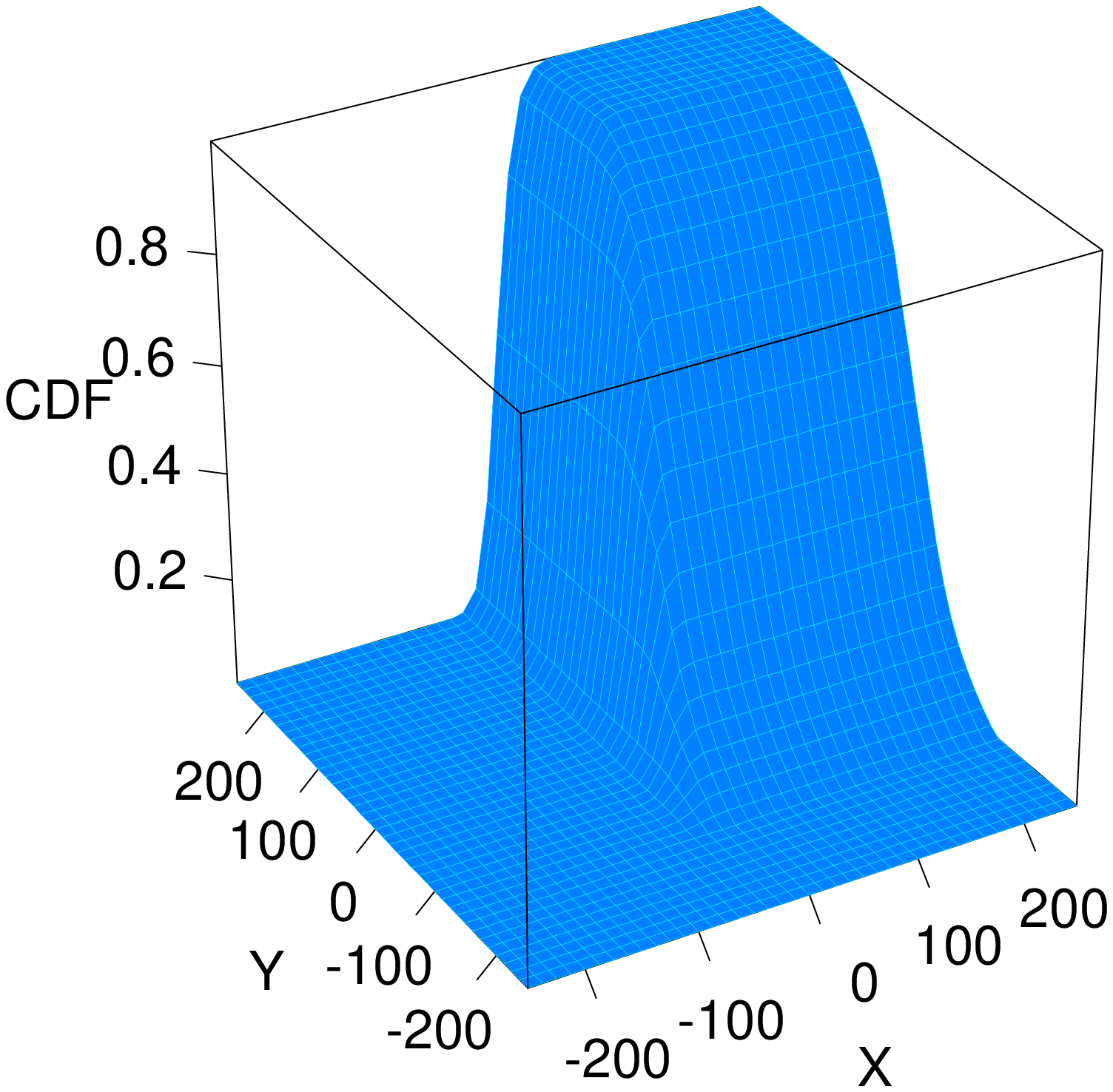,width=\textwidth}
\end{minipage}
\hfill
\begin{minipage}{4cm}
\epsfig{figure=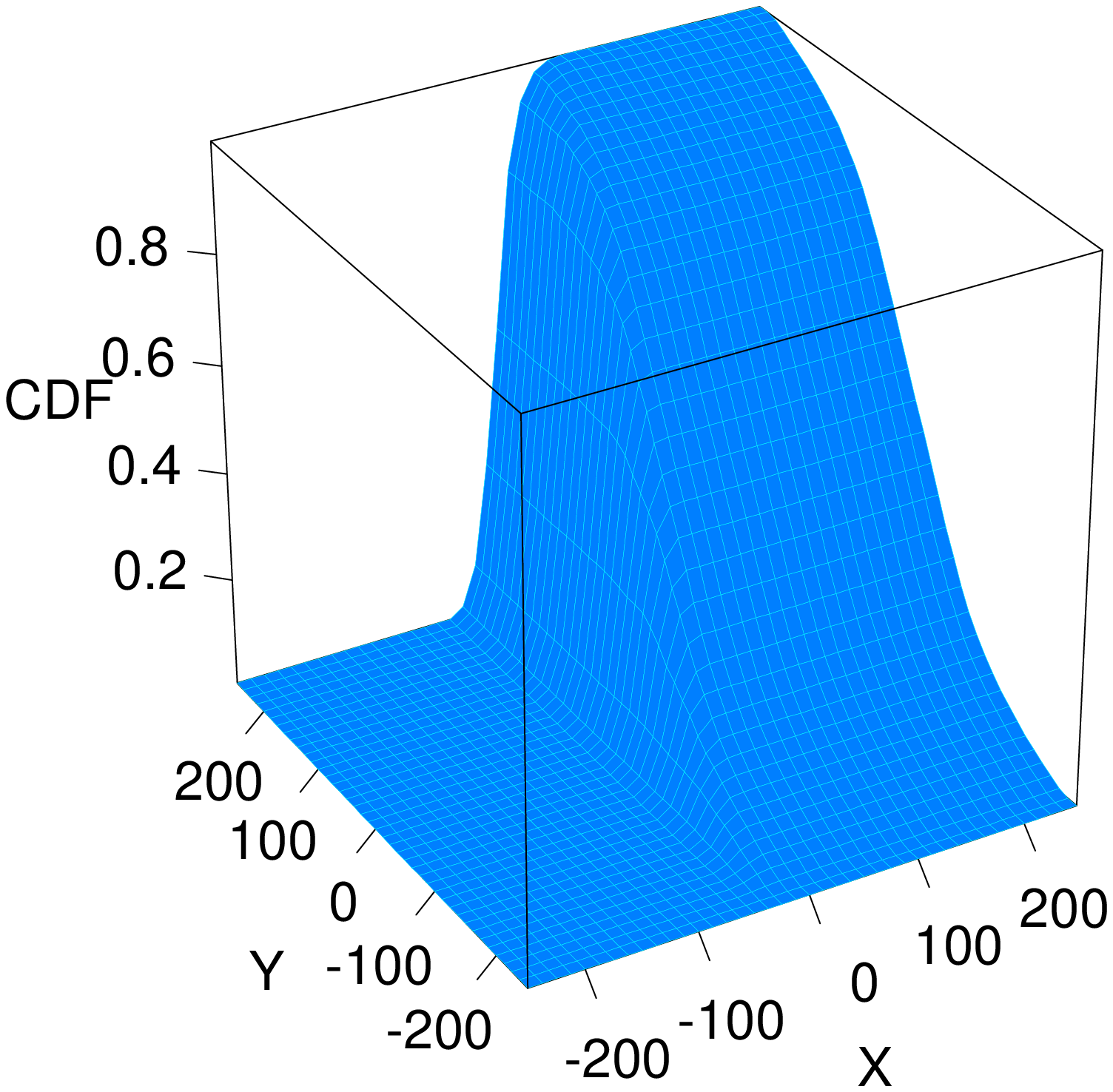,width=\textwidth}
\end{minipage}
\hfill
\caption{CDFs at different time instances computed by the CPI method; 
top left: $t=0$, top right: $t=300 \Delta t$, bottom left: $t=600 \Delta t$, 
bottom right: $t=900 \Delta t$.} 
\label{fig5:fig}
\end{figure}

\begin{figure}
\centerline{\epsfig{figure=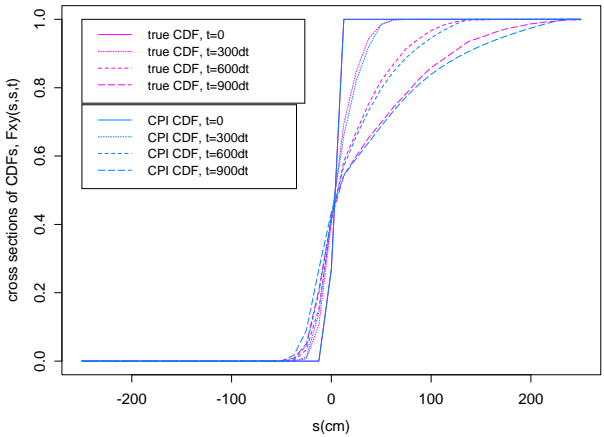,width=0.35\textwidth,angle=-90}}
\caption{Comparison between cross sections of true CDFs 
and CDFs computed by CPI.}
\label{fig6:fig}
\end{figure}

We now use CPI to accelerate the evolution of the coarse-scale observables. 
The coarse-scale initial condition is the same as that in Simulation 1.
Particles are evolved for an initial block of 10 ($l=10$) 
fine-scale time steps and then again for another block of 10 ($n-l=10$) 
fine-scale time steps. 
At each of the latter 10 steps, the marginal ICDF and 
20 ($M=20$) conditional ICDFs of particle positions are formed; 
the time series of the coefficients of their leading
modes is linearly extrapolated (with a slope estimated through least-squares) 
over a time interval equal to 10 ($T=10$) fine-scale time steps. 
The basis onto which the ICDFs are projected consists here of
shifted Legendre polynomials of order up to and including 5 ($P=5$). 
Since the ICDFs are anti-symmetric with respect to the axis $f=0.5$, 
the coefficients of the 2nd and 4th modes vanish.
Hence, only 84 coefficients need to be extrapolated. 
At the end of the extrapolation, the ICDFs are reconstructed 
and particle positions, generated according to these ICDFs, are simulated again. 
At time $300\Delta t$, $600 \Delta t$ and $900 \Delta t$, particle 
positions are recorded and CDFs plotted in Fig. \ref{fig5:fig}. 
The number of grid points used to compute CDFs is again 1681. 
The cross sections, $F_{XY}(s,s,t)$, of CDFs in Simulation 1 
and 2 are compared in Fig. \ref{fig6:fig}, which shows  
an excellent visual match between the true CDFs and those computed from the CPI algorithm.

\bigskip
{\noindent\bf Simulation 3: Coarse Dynamic Renormalization of the 
Self-Similar Particle System}

Pretending that the macroscale equation (\ref{eqn12:eqn}) is not available, 
we now directly use the microscale simulator (\ref{eqn10:eqn}) 
to compute the constants $p$ and $a$, the similarity exponent $\alpha$ 
and the macroscopic self-similar solution.

In the approach provided in Section 4, 
Newton's method is used to solve the equation (\ref{eqnpsolve:eqn}) for $p$. 
The test function $f$ was first chosen as a 2-dimensional 
joint Gaussian distribution function, $f(x,y)=1/16N(x/4)N(y/4)$, 
where $N(x)$ and $N(y)$ are standard Gaussian distributions. 
We select the value of the positive real number $A=2.0$. 
The two coordinates $(u_1,v_1)$ and $(u_2,v_2)$ in (\ref{eqnpsolve:eqn}) 
are chosen as $(-2,-2)$ and $(3,3)$, respectively. 
To reduce fluctuations of values for the operator 
$D_{xy}$, $9000$ particles are used and $500$ replica
copies of values for $D_{xy}$ and $D_{uv}$ are averaged in the computation.

Starting from the initial value $p_0=5.0$, 
iterative values for $p$ are stabilized at 3.0 after 3 iterations. 
Accordingly, the converged value for $a$ is $-2.0$ (Table \ref{tabpa1:tab}). 

\begin{table}
\begin{center}
\begin{tabular*}{0.45\textwidth}%
   {@{\extracolsep{\fill}}cccr}
   No. of Iterations &    $p$   &    $a$ \cr
   \hline
   0                  &   5.0  & -3.34959  \cr
   1                  &  2.80093  & -1.88311  \cr
   2                  &  2.99246  & -2.03967  \cr
   3                  &  3.00106  & -2.04212  \cr
   4                  &  3.00370  & -2.05468  \cr
   5                  &  2.99592  & -2.03927  \cr
   6                  &  2.99659  & -2.03914  \cr
   7                  &  2.99753  & -2.04017  \cr
   8                  &  2.99831  & -2.04189  \cr
   \hline
\end{tabular*}
\end{center}
\caption{Iterative evaluation of the constants $p$ and $a$.}
\label{tabpa1:tab}
\end{table}

For verification, we choose another set of parameters 
to determine $p$ and $a$: $f(x,y)=1/25N(x/5)N(y/5)$, 
$A=2.5$, $(u_1,v_1)=(-3,-3)$, and $(u_2,v_2)=(4,4)$. 
Again, starting from the initial value $p_0=5.0$, 
iterative evaluations of $p$ are stabilize at $3.0$ 
after 3 iterations and the converged value for $a$ 
is again around $-2.0$ (Table \ref{tabpa2:tab}).

\begin{table}
\begin{center}
\begin{tabular*}{0.45\textwidth}%
   {@{\extracolsep{\fill}}cccr}
   No. of Iterations &    $p$   &    $a$ \cr
   \hline
   0                  &   5.0  & -3.38239  \cr
   1                  &  3.19621  & -2.25167  \cr
   2                  &  2.97334  & -2.06146  \cr
   3                  &  2.99987  & -2.08557  \cr
   4                  &  2.99574  & -2.08444  \cr
   5                  &  2.99347  & -2.08425  \cr
   6                  &  2.99747  & -2.08345  \cr
   7                  &  2.99247  & -2.08073  \cr
   8                  &  2.99833  & -2.08601  \cr
   \hline
\end{tabular*}
\end{center}
\caption{Iterative values for the constants $p$ and $a$ for a different set of 
algorithm parameters (see text).}
\label{tabpa2:tab}
\end{table}
 
We can therefore conclude that there exist $p=3.0$ and $a=-2.0$ 
such that the unavailable differential operator $D_{xy}$ 
corresponding to the microsimulator (\ref{eqn9:eqn}) 
possesses the scale invariance property (\ref{eqnDxy:eqn}). 
Accordingly, the similarity exponent is $\alpha=0.5$ by Eqn. (\ref{eqnalphaa:eqn}).

The template condition for the $x$ direction is chosen to be $\omega(-2.832,\infty,t)=0.4$, 
i.e., the u-coordinate corresponding to the renormalized marginal 
CDF $\omega_U=0.4$ always has the same value, $-2.832$cm. 
The constant $c$ in the analytical solution (\ref{eqn104:eqn}) 
is obtained as $c=0.2sec^{-1}$ based on our template. 
The corresponding std.'s for the analytical 
self-similar shape are $\sigma_X=5\sqrt{5}$cm 
and $\sigma_Y=25\sqrt{15}/3$cm, respectively.  

\begin{figure}
\begin{minipage}{4cm}
\epsfig{figure=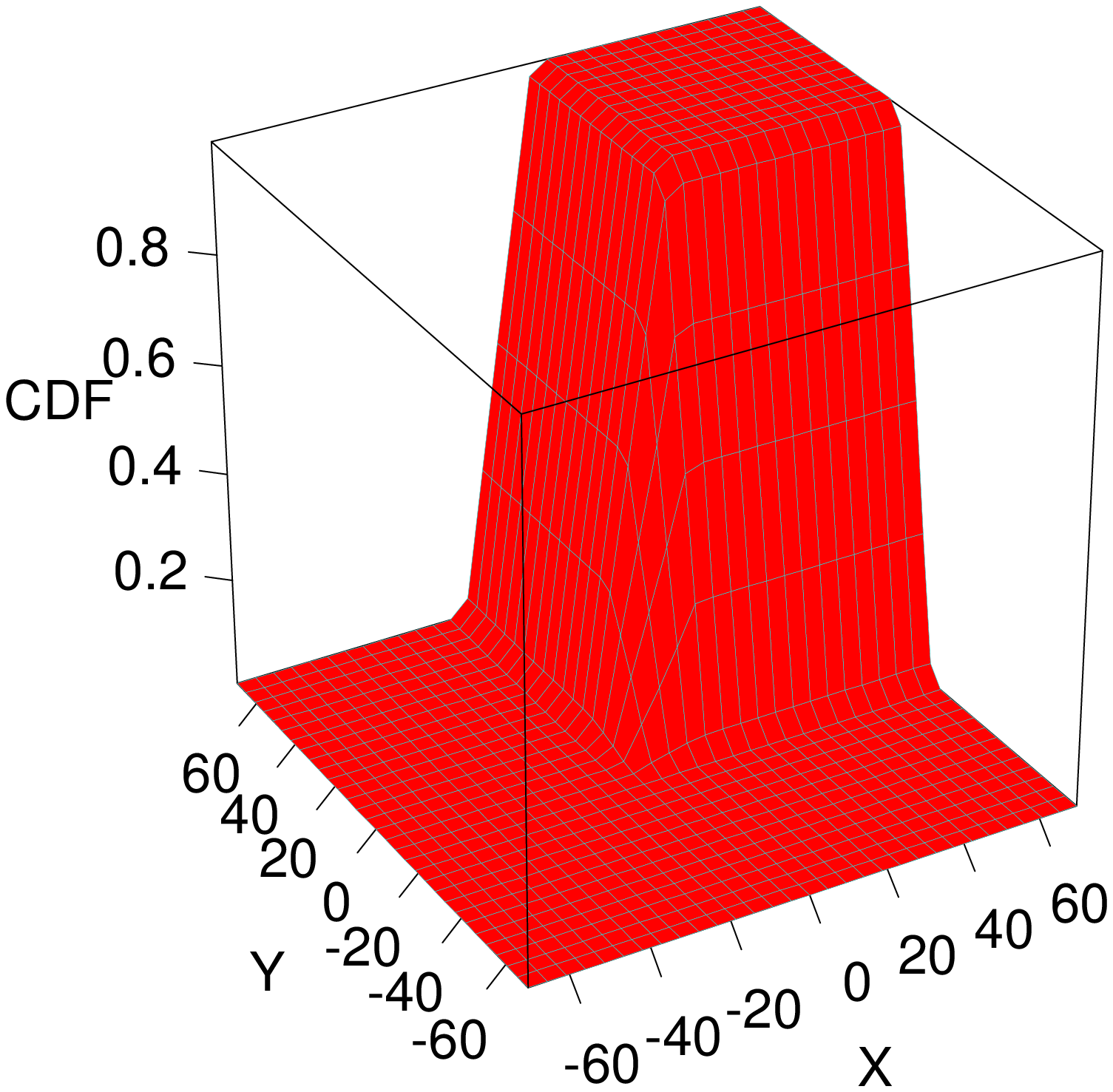,width=\textwidth}
\end{minipage}
\hfill
\begin{minipage}{4cm}
\epsfig{figure=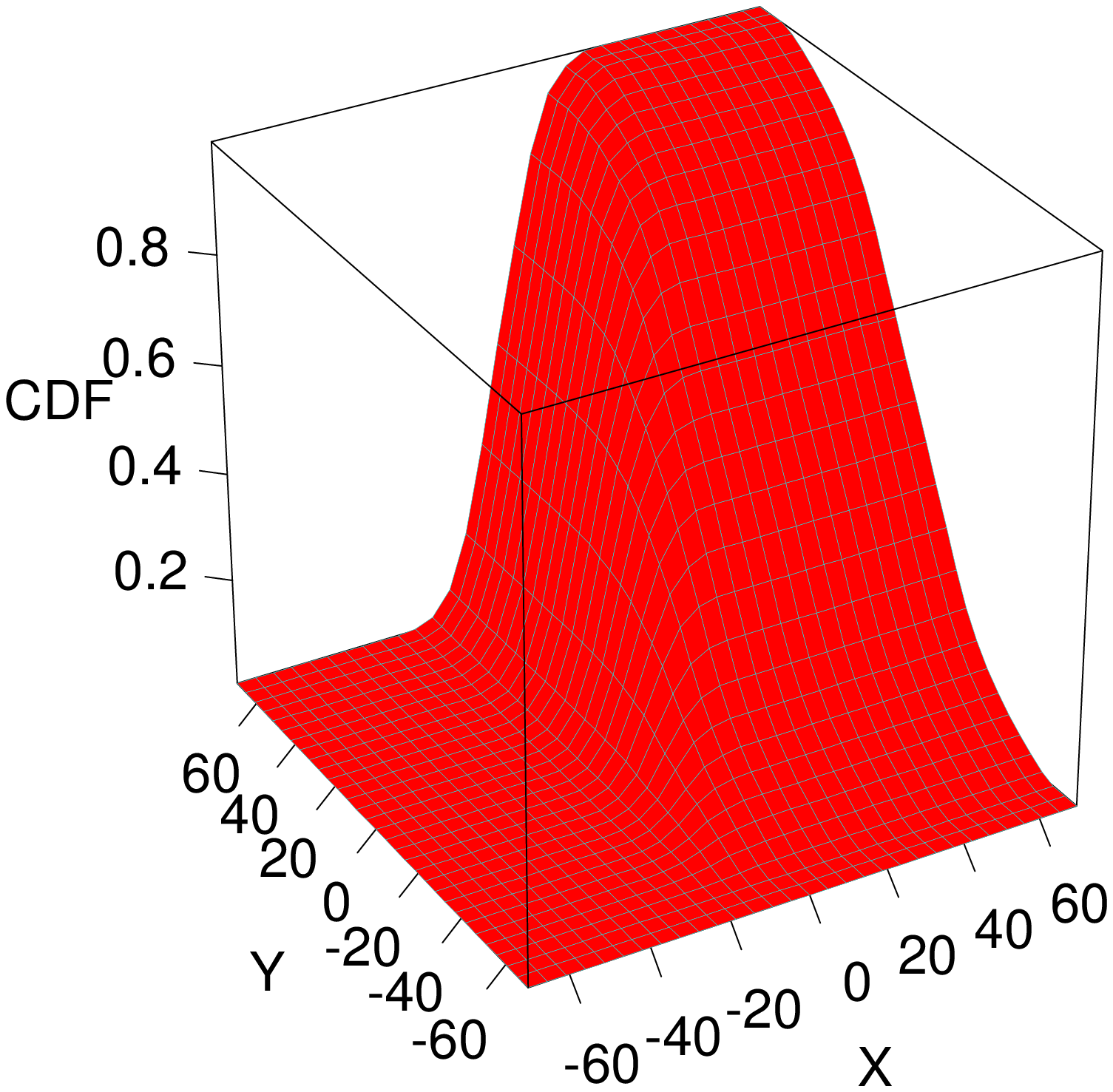,width=\textwidth}
\end{minipage}
\hfill
\begin{minipage}{4cm}
\epsfig{figure=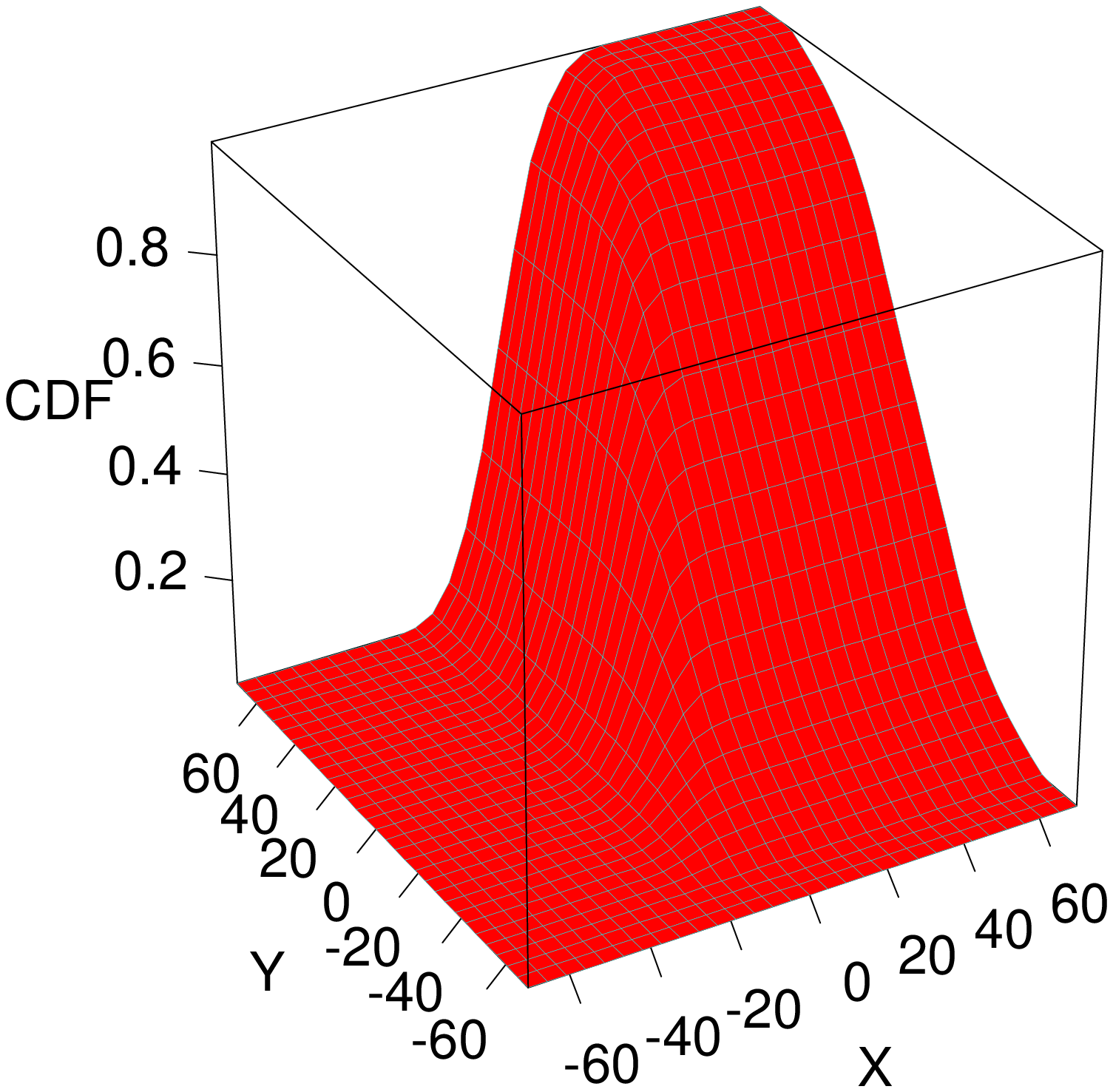,width=\textwidth}
\end{minipage}
\hfill
\begin{minipage}{4cm}
\epsfig{figure=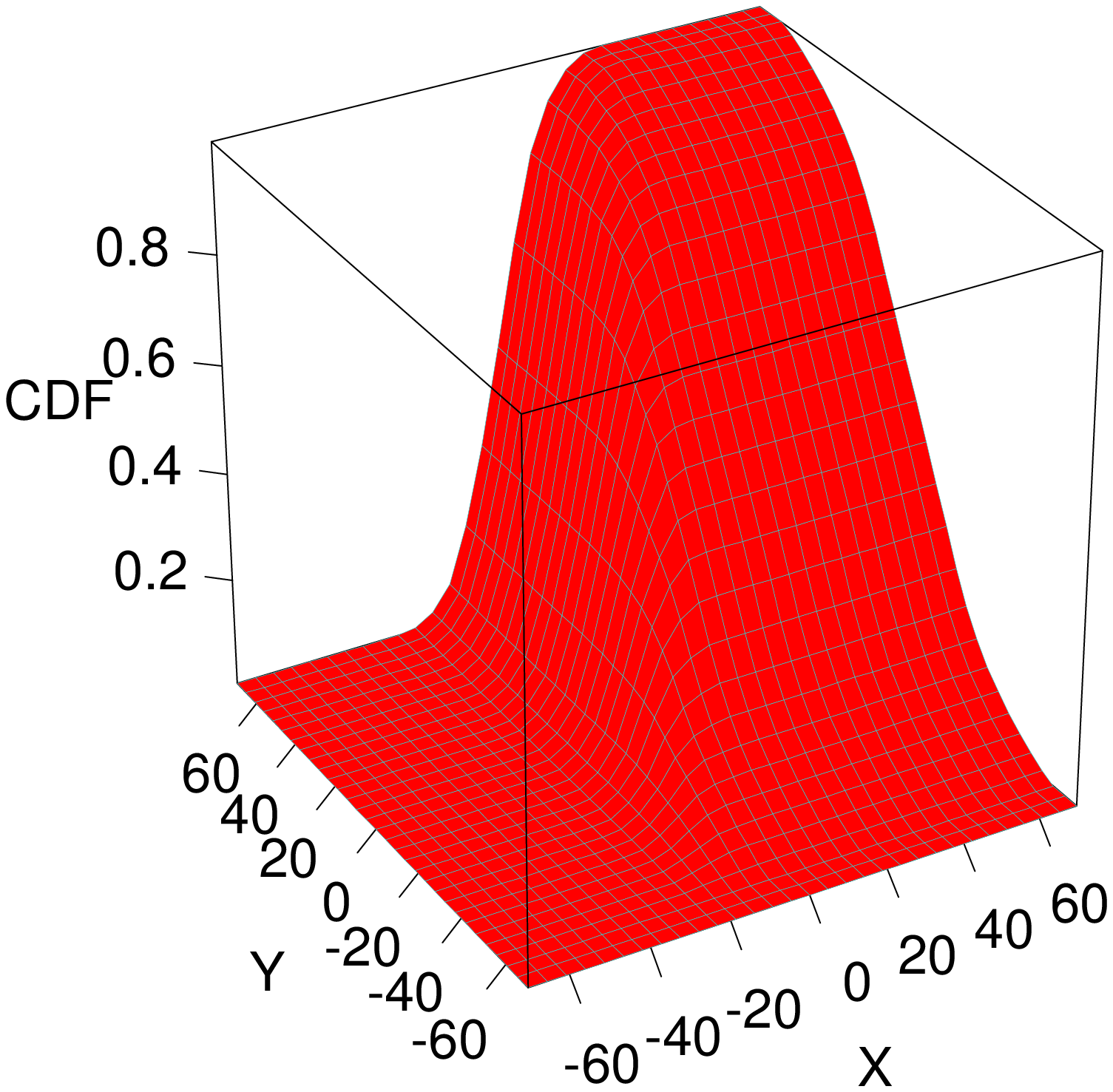,width=\textwidth}
\end{minipage}
\hfill
\caption{Renormalized CDFs in Simulation 3; 
top left: initial CDF, top right: renormalized CDF 
after 2nd iteration, bottom left: renormalized CDF 
after 4th iteration, bottom right: renormalized CDF after 6th iteration.} 
\label{fig8:fig}
\end{figure}

\begin{figure}
\centerline{\epsfig{figure=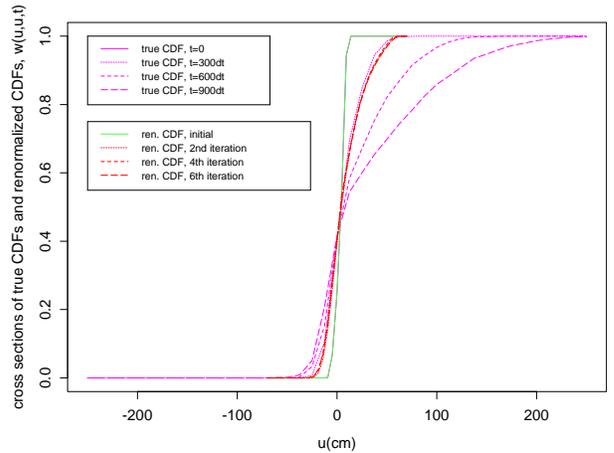,width=0.35\textwidth,angle=-90}}
\caption{Comparison between cross sections of true CDFs and 
rescaled CDFs in Simulation 3.}
\label{fig9:fig}
\end{figure}

The CDF corresponding to a uniform distribution of particle positions 
over the space domain $(-10cm,10cm)\times(-10cm,10cm)$ is 
used as the initial condition. 
Direct iteration is used to solve for the fixed point of equation (\ref{eqn8:eqn}). 
The time interval $T'$ is $100 \Delta t$. 
The number of conditional ICDFs is 20 ($M=20$) and 
the basis for the ICDFs is again the shifted Legendre polynomials 
of order up to and including 5 ($P=5$). 
In this simulation, 100 copies of ensemble particle 
positions are generated according to the mode coefficients of the ICDFs 
at the beginning of each iteration and let to evolve.
The mode coefficients at the end of each iteration are 
obtained by averaging over these 100 replica copies. 
After the 2nd, 4th and 6th iterations, renormalized mode 
coefficients of the ICDFs are used to generate particle positions, 
out of which the CDFs are computed and plotted respectively in Fig. \ref{fig8:fig}. 
Fig. \ref{fig9:fig} also compares the cross sections $\omega(u,u,t)$ 
of true CDFs and renormalized CDFs in this simulation. 
Clearly, the renormalized solutions quickly approach the self-similar steady state. 

\begin{figure}
\centerline{\epsfig{figure=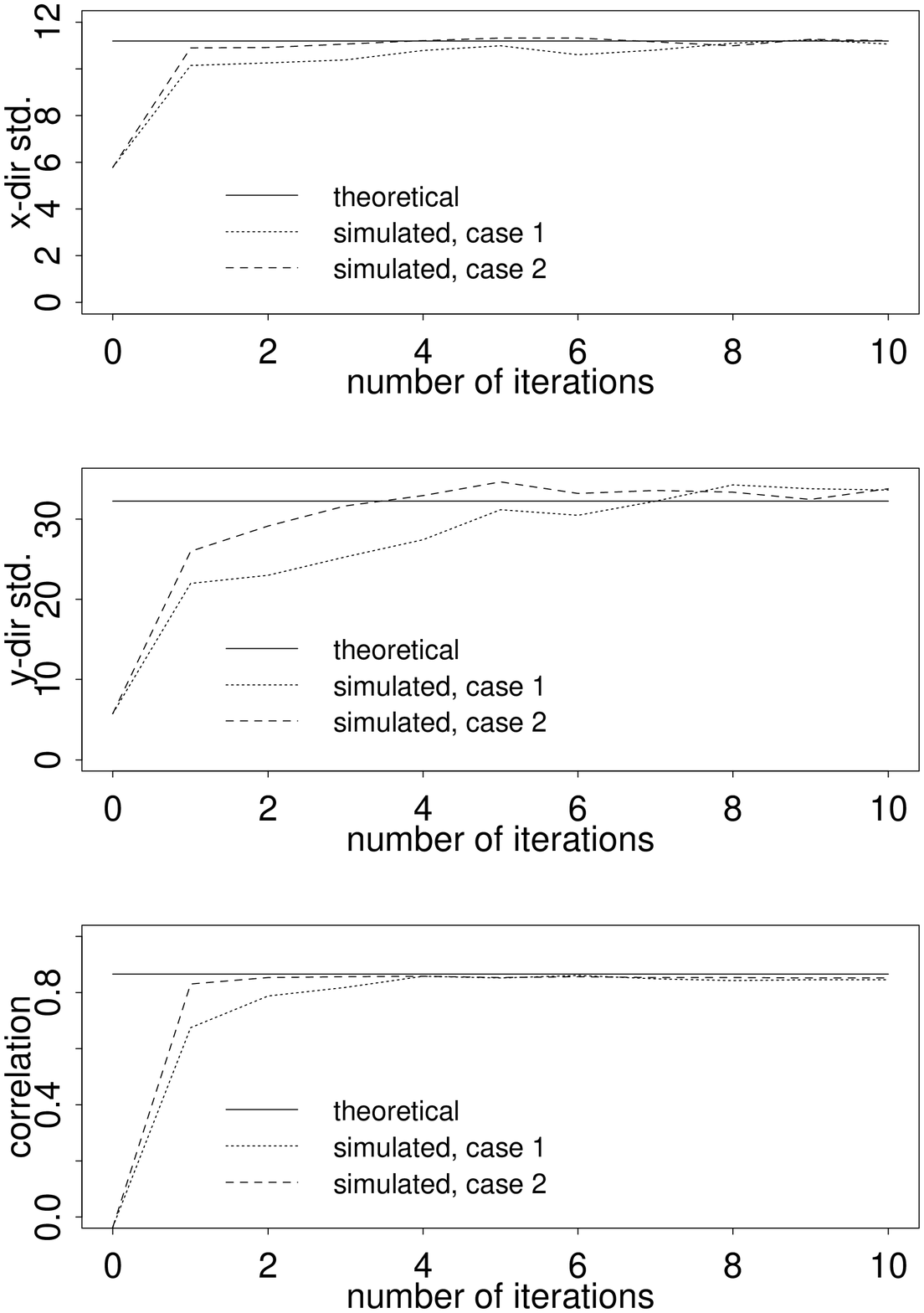,width=0.4\textwidth}}
\caption{Comparison for standard deviations 
and correlations of self-similar shapes in Simulation 3: Cases 1,2.}
\label{fig10:fig}
\end{figure}

To validate the computation of the self-similar solution shape, the std.'s and 
correlation of the computed shape are compared with those of 
the known analytical solution.
The std.'s and correlations of the rescaled CDFs are calculated 
via the ensemble particle positions corresponding to these CDFs. 
The comparison is shown in Fig. \ref{fig10:fig}, where curves 
in Case 1 represent results obtained using this template condition and time interval.
The std.'s and correlations of the rescaled CDFs approach those of the 
analytical self-similar shape, which means that the rescaled CDF 
coincides eventually with a member in the family of theoretical 
self-similar shapes expressed by Equation (\ref{eqn104:eqn}).
 
As the renormalized CDF $\omega$ reaches its steady state, 
we can set this CDF as the initial condition and evolve the 
microscale dynamics (\ref{eqn10:eqn}) for two more loops with $t_1=100\Delta t$ 
and $t_2=300 \Delta t$. 
The rescaling variable $A(t)$ is listed in Table \ref{tabA:tab}. 
Note that $A(t)=1$ at $t=0$. 
By Equation (\ref{alphacal:eqn}), the similarity exponent $\alpha$ is approximated
as $0.520$, within $4 \%$ of the theoretical value $1/2$.

\begin{table}
\begin{center}
\begin{tabular*}{0.45\textwidth}%
   {@{\extracolsep{\fill}}cccr}
   $t(sec)$           &   $A(t)$  &  $A_t(t)$ \cr
   \hline
   0                  &  1.00000  &    -      \cr
   1                  &  1.10268  &  0.10268  \cr
   3                  &  1.27793  &  0.08763  \cr
   \hline
\end{tabular*}
\end{center}
\caption{The rescaling variable $A(t)$.}
\label{tabA:tab}
\end{table}

In the following, the effect of variation of templates and 
evolution times of the fixed-point operator $\Phi_{T'}$, 
on the computed renormalized self-similar shapes will be examined. 
For self-similar systems, we can see from Eqn. (\ref{eqnomega:eqn}) 
that, as the system reaches the steady state, the 
rescaled shape of CDFs will remain the same, irrespective of changes in
the evolution time $T'$. 
Also, the steady-state CDF shapes will coincide with members 
in the family of self-similar solutions prescribed by 
Eqn. (\ref{eqnselfsimilar:eqn}) no matter what the scale of the template is. 

We choose four cases of the template condition and evolution 
time including the one above:

\begin{enumerate}
\item $\omega(-2.832,\infty,t)=0.4$, $T'=100\Delta t$;
\item $\omega(-2.832,\infty,t)=0.4$, $T'=200\Delta t$;
\item $\omega(-0.283,\infty,t)=0.4$, $T'=100\Delta t$;
\item $\omega(-0.283,\infty,t)=0.4$, $T'=200\Delta t$.
\end{enumerate} 

\begin{figure}
\centerline{\epsfig{figure=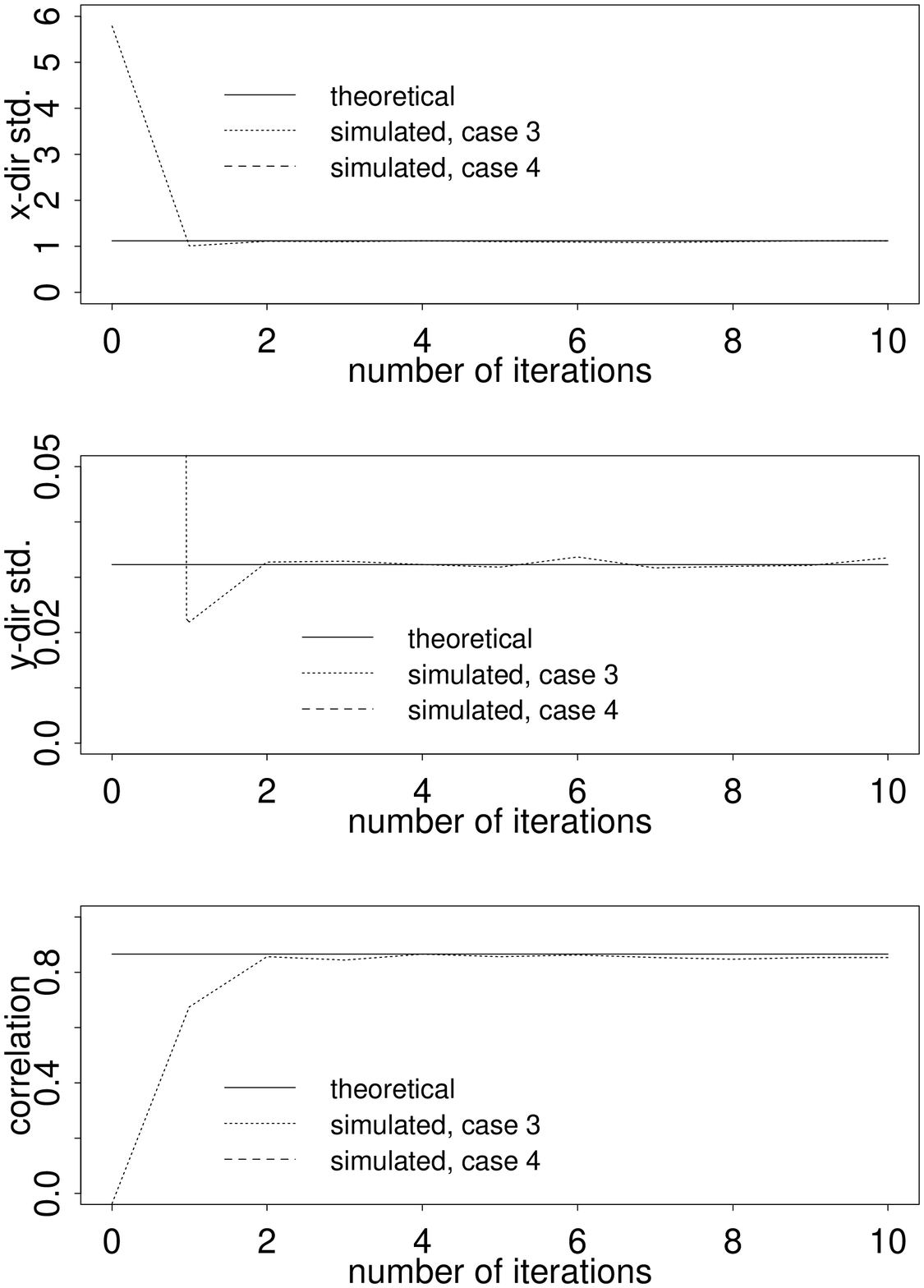,width=0.4\textwidth}}
\caption{Comparison for standard deviations and 
correlations of self-similar shapes in Simulation 3: Cases 3,4.}
\label{selfsimtemevocom2:fig}
\end{figure}

The iterative values of std.'s and correlation for the 
four cases are shown in figures \ref{fig10:fig} 
and \ref{selfsimtemevocom2:fig}. 
Comparison with theoretical calculations shows 
that variation of templates and evolution times 
indeed does not cause deviation of the converged rescaled CDF 
from the family of self-similar solutions.    
     
\bigskip
{\noindent\bf Simulation 4: Coarse Dynamic Renormalization of the 
Asymptotically Self-Similar Particle Dynamics}

For the particle system in Section 5.2, 
the procedure in Section 4 is used to 
check if its macroscopic differential operator $D_{xy}$ 
possesses the scale invariance property (\ref{eqnDxy:eqn}). 
The two parameter sets in Simulation 3 are used here.  

\begin{enumerate}
\item Parameter set 1: $f(x,y)=1/16N(x/4)N(y/4)$, $A=2.0$, $(u_1,v_1)=(-2,-2)$,  $(u_2,v_2)=(3,3)$;
\item Parameter set 2: $f(x,y)=1/25N(x/5)N(y/5)$, $A=2.5$, $(u_1,v_1)=(-3,-3)$,  $(u_2,v_2)=(4,4)$.
\end{enumerate}

Newton's method is utilized again to solve for $p$ and $a$ 
for each parameter set. 
The iteratively computed values of $p$ and $a$ are 
listed in tables \ref{tabpa4:tab} and \ref{tabpa5:tab}.

\begin{table}
\begin{center}
\begin{tabular*}{0.45\textwidth}%
   {@{\extracolsep{\fill}}cccr}
   No. of Iterations &    $p$   &    $a$ \cr
   \hline
   0                  &   5.0  & -3.87764  \cr
   1                  &  4.14589  & -3.32940  \cr
   2                  &  3.81151  & -3.09345  \cr
   3                  &  3.78553  & -3.08460  \cr
   4                  &  3.79408  & -3.08433  \cr
   5                  &  3.78576  & -3.08486  \cr
   6                  &  3.79758  & -3.08472  \cr
   7                  &  3.78701  & -3.09319  \cr
   8                  &  3.79517  & -3.08471  \cr
   \hline
\end{tabular*}
\end{center}
\caption{Iteratively computed values of $p$ and $a$ for Parameter set 1
(see text).}
\label{tabpa4:tab}
\end{table}

\begin{table}
\begin{center}
\begin{tabular*}{0.45\textwidth}%
   {@{\extracolsep{\fill}}cccr}
   No. of Iterations &    $p$   &    $a$ \cr
   \hline
   0                  &   5.0  & -4.23757  \cr
   1                  &  4.47814  & -3.85506  \cr
   2                  &  2.58903  & -1.95833  \cr
   3                  &  4.01019  & -3.38693  \cr
   4                  &  3.64906  & -2.96049  \cr
   5                  &  3.35744  & -2.69786  \cr
   6                  &  3.39829  & -2.74037  \cr
   7                  &  3.40776  & -2.75770  \cr
   8                  &  3.41325  & -2.74170  \cr
   9                  &  3.41298  & -2.75793  \cr
   \hline
\end{tabular*}
\end{center}
\caption{Iteratively computed values of $p$ and $a$ for Parameter set 2
(see text).}
\label{tabpa5:tab}
\end{table}

It can be seen that now converged values of $p$ and $a$ {\it do vary} 
with the template scale and value of $A$. 
As the template scale and $A$ increase, values of $p$ and $a$ 
approach those in the self-similar case. 
Pretending that the macroscale equation is not explicitly known, 
we may suspect that the particle system exhibits asymptotically self-similar
dynamics. 
Therefore, for ``asymptotically large enough" template conditions, 
the operator $D_{xy}$ still approximately possesses the scale invariance 
property (\ref{eqnDxy:eqn}); as we did in the self-similar case, 
we can use a fixed point algorithm to find a 
long-time steady state for the asymptotically self-similar solution. 
In analogy to the self-similar case, 
the evolution time interval in the fixed-point operator 
does not affect the converged shape for such large enough scales.
%
%
%
%

\begin{figure}
\centerline{\epsfig{figure=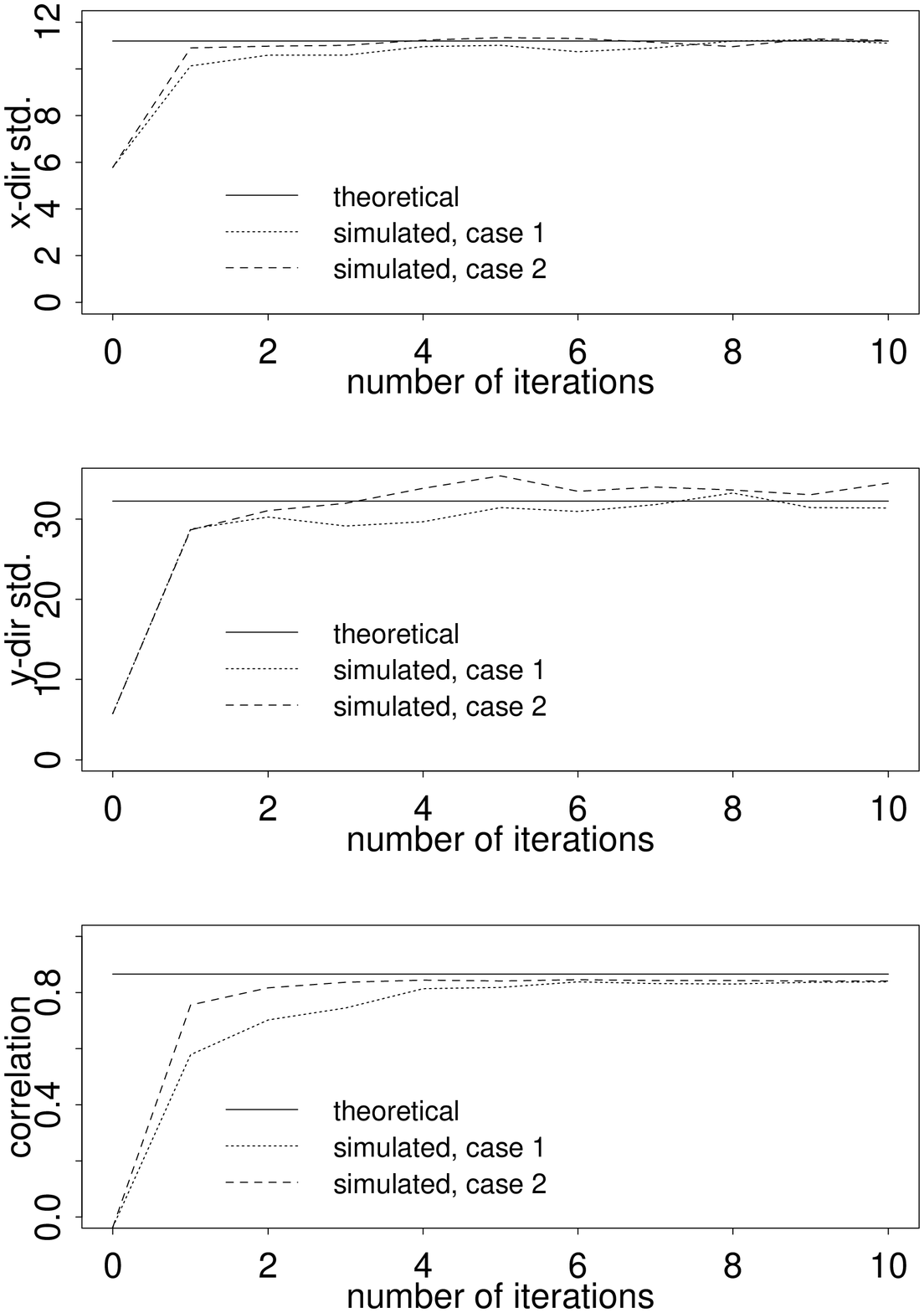,width=0.4\textwidth}}
\caption{Comparison for standard deviations and correlations of 
converged shapes in Simulation 4: Cases 1,2.}
\label{asysimvarcor1:fig}
\end{figure}

\begin{figure}
\centerline{\epsfig{figure=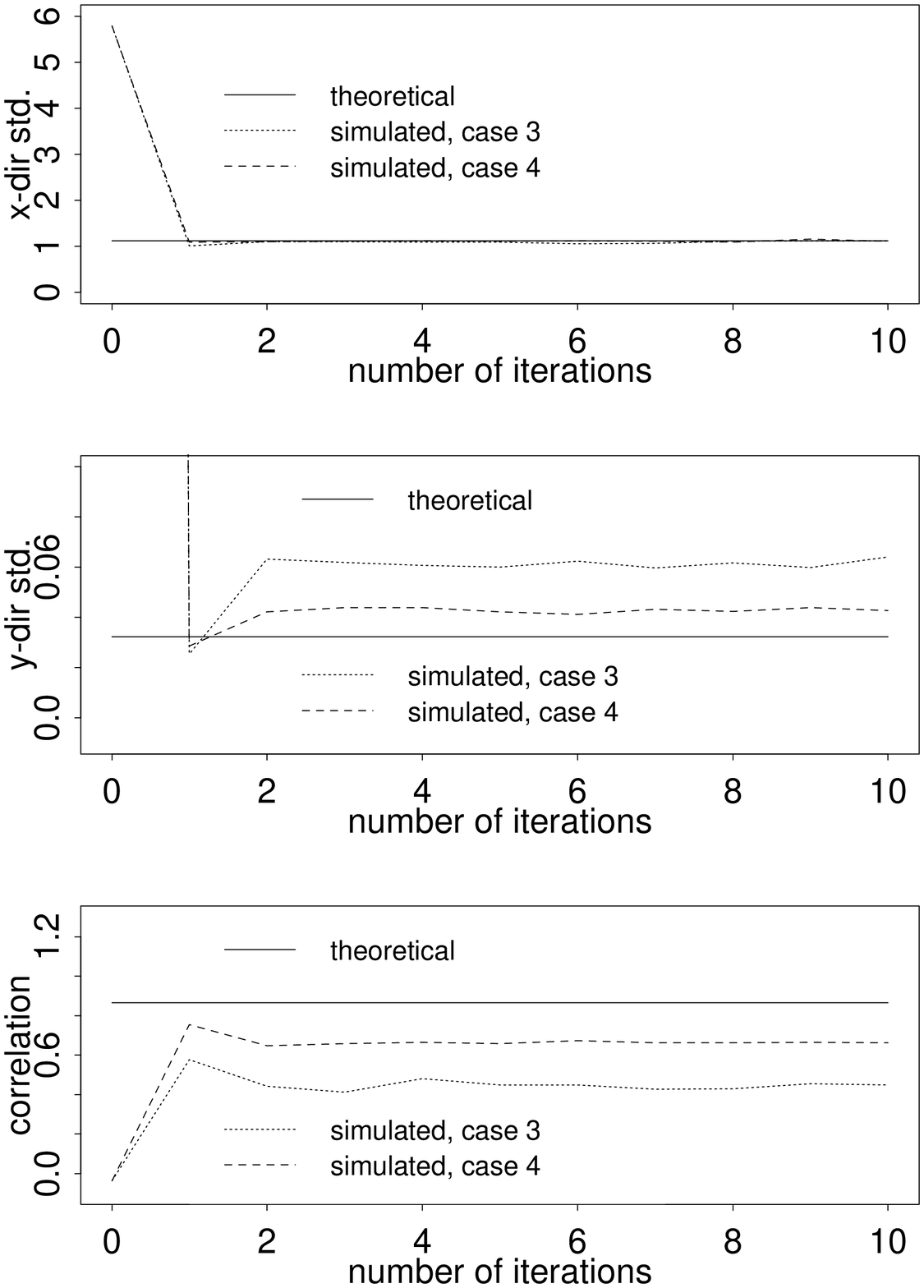,width=0.4\textwidth}}
\caption{Comparison for standard deviations and correlations of 
converged shapes in Simulation 4: Cases 3,4.}
\label{asysimvarcor2:fig}
\end{figure}

The four template condition and evolution time cases in Simulation 3 
are used in the fixed point algorithm to verify the above assertions. 
The value $p$ is set to 3.0,
the same value as that in Simulation 3. 
The microsimulator is the discretized dynamics 
(\ref{dis2ddifBrownian:eqn}). 
Iterative values of std.'s and correlation 
for the four cases are shown in figures \ref{asysimvarcor1:fig} and 
\ref{asysimvarcor2:fig}. 
As can be seen, for large template conditions, 
the length of evolution time $T'$ does not 
affect the converged values of std.'s and 
correlation for distribution of particle positions. 
Yet for small templates, the effect of $T'$ on the 
$y$-direction std. and correlation is evident. 
As $T'$ increases from $100\Delta t$ to $200\Delta t$, 
the converged values of $y$-direction std. and 
correlation approach their expected theoretical values.
%
%

\vspace{0.5cm}
\noindent{\bf 7. Conclusions and Remarks} 
\vspace{0.2cm}

We presented an equation-free computational approach, 
based on using marginal and conditional ICDFs as 
coarse-scale observables, for the computer-assisted study of 
multidimensional random particle system dynamics. 
Coarse projective integration employing this time-stepper can be applied to accelerate the
computational evolution of particle CDF computations; the approach targets 
multidimensional particle systems whose coarse-scale models are not explicitly available.
%
%
Coarse dynamic renormalization can also be used to 
analyze particle systems with self-similar or asymptotically self-similar 
coarse-grained evolution dynamics, and to
obtain long-time renormalized steady state (self-similar) solutions.

The examples in this paper are admittedly simple, 
yet they illustrate the computational approaches in a context where the
results can be validated; we hope that the type of
multiscale algorithms presented here may be useful in more 
complicated situations
(e.g., particles mixing in time-dependent velocity fields)
if the macroscopical dynamic are effectively self-similar.
Another possible application of such equation-free approaches is in
cases where even {\it the coarse-scale observables} are 
characterized by uncertainty/stochasticity. 
Polynomial chaos observables 
have been used in the solution of explicit macroscale PDEs 
for passive scalar transport, where the uncertainty enters
through random initial conditions or 
boundary conditions \cite{Ghanem:98,Xiu:03}.
Such polynomial chaos observables may be combined with
the coarse-graining techniques presented here when no
explicit coarse-scale descriptions of the particle 
system dynamics are available.

\vspace{0.5cm}
\noindent{\bf Appendix A}
\vspace{0.2cm}

For Equation (\ref{eqn11:eqn}) with the initial condition
\begin{equation}
   P_{XY}(x,y,t_0)=\delta(x)\delta(y),  \nonumber
\end{equation}
(where $t_0$ is the blowup time) it is shown in the following that the solution is
\begin{equation}
P_{XY}(x,y,t)={\sqrt{3} \over {\pi {D^2} (t-t_0)^2 }} e^{- \left( {{6(y - 0.5x(t-t_0))^2} \over {D^2 (t-t_0)^3}} + {{x^2} \over {2D^2 (t-t_0)}} \right)}. \nonumber 
\end{equation}

We can see that
\begin{eqnarray}
 &&{{\partial P_{XY}} \over {\partial t}} = e^{- \left( {{6(y - 0.5x(t-t_0))^2} \over {D^2 (t-t_0)^3}} + {{x^2} \over {2D^2(t-t_0)}} \right)} \nonumber
  \\
 && \cdot[{{-2\sqrt{3}} \over {\pi D^2 (t-t_0)^3}} + {{\sqrt{3}} \over {\pi D^2 (t-t_0)^2}} ( {{(x-y)^2} \over {2D^2 (t-t_0)^2}} \nonumber
  \\
 && +{{18(y-0.5x(t-t_0))^2} \over {D^2(t-t_0)^4}} +{{6(y-0.5x(t-t_0))x} \over {D^2 (t-t_0)^3}} ) ], \nonumber
\end{eqnarray}

\begin{eqnarray}
 &&x{{\partial P_{XY}} \over {\partial y}} = e^{- \left( {{6(y - 0.5x(t-t_0))^2} \over {D^2 (t-t_0)^3}} + {{x^2} \over {2D^2(t-t_0)}} \right)} \nonumber
                     \\
 &&\cdot \left[{{-\sqrt{3}x} \over {\pi D^2 (t-t_0)^2}} \cdot {{12(y-0.5x(t-t_0))} \over {D^2 (t-t_0)^3}} \right], \nonumber
\end{eqnarray}

\begin{eqnarray}                                
 &&{{D^2} \over 2} {{\partial^2 P_{XY}} \over {\partial x^2}} = {{D^2} \over 2} e^{- \left( {{6(y - 0.5x(t-t_0))^2} \over {D^2 (t-t_0)^3}} + {{x^2} \over {2D^2(t-t_0)}} \right)} \nonumber
                             \\
 &&\cdot \left[ {{-4\sqrt{3}} \over {\pi D^2 (t-t_0)^3}} + {{\sqrt{3}} \over {\pi D^4 (t-t_0)^4}} \cdot{{(6(y)-4x(t-t_0))^2} \over {D^2(t-t_0)^2}} \right]. \nonumber
\end{eqnarray}   

Hence, 
\begin{equation}
  {{\partial P_{XY}} \over {\partial t}} + x {{\partial P_{XY}} \over {\partial y}} - {{D^2}\over 2} {{\partial^2 P_{XY}} \over {\partial x^2}}= 0.  \nonumber
\end{equation}

\clearpage


\end{document}